\documentclass[nonblindrev,copyedit]{informs3}
\usepackage{setspace,graphics,booktabs,color,enumerate,subfigure,url,bm}
\usepackage{appendix}
\usepackage{amssymb}
\usepackage{mathrsfs}
\usepackage{multirow}
\usepackage{amsmath, amssymb, color}
\usepackage{graphicx}
\usepackage{subfigure}
\usepackage{slashbox}
\usepackage{algorithmic}
\definecolor{DarkBlue}{rgb}{0,0.08,0.45}
\usepackage[backref = false, bookmarks, colorlinks = true, plainpages = false, citecolor = DarkBlue , urlcolor = DarkBlue, filecolor = DarkBlue, linkcolor = DarkBlue]{hyperref}
\usepackage{array}
\usepackage{float}
\floatstyle{ruled}
\newfloat{algorithm}{htp}{lop}
\floatname{algorithm}{Algorithm}

\newtheorem{example}{Example}

\newcommand{\cV}{\mbox{\boldmath $\mathcal{V}$}}
\newcommand{\OPT}{{\mbox{OPT}}}

\newcommand{\bx}{\mbox{\boldmath $x$}}
\newcommand{\by}{\mbox{\boldmath $y$}}
\newcommand{\bz}{\mbox{\boldmath $z$}}

\usepackage{natbib}
 \bibpunct[, ]{(}{)}{,}{a}{}{,}%
 \def\bibfont{\small}%

\usepackage{tikz}   
\usepackage{mathdots}
\usetikzlibrary{arrows}
\usetikzlibrary{decorations.pathreplacing}

\theoremstyle{TH}%
\newtheorem{theorem}{Theorem}
\newtheorem{lemma}{Lemma}

\newtheorem{corollary}{Corollary}

\ECRepeatTheorems

\EquationsNumberedThrough    

\MANUSCRIPTNO{}

\begin{document}




\TITLE{\large{Assortment Optimization under a Single Transition
Model}}

\ARTICLEAUTHORS{\AUTHOR {Kameng Nip} \AFF{Department of Mathematical
Sciences, Tsinghua University, Beijing, China
\href{mailto:njm13@mails.tsinghua.edu.cn}{njm13@mails.tsinghua.edu.cn}}
\AUTHOR {Zhenbo Wang} \AFF{Department of Mathematical Sciences,
Tsinghua University, Beijing, China
\href{mailto:zwang@math.tsinghua.edu.cn}{zwang@math.tsinghua.edu.cn}}
\AUTHOR {Zizhuo Wang} \AFF{Department of Industrial and Systems
Engineering, University of Minnesota, MN, USA
\href{mailto:zwang@umn.edu}{zwang@umn.edu}}}

\date{\today}

\ABSTRACT{In this paper, we consider a Markov chain choice model
with single transition. In this model, customers arrive at each
product with a certain probability. If the arrived product is
unavailable, then the seller can recommend a subset of available
products to the customer and the customer will purchase one of the
recommended products or choose not to purchase with certain
transition probabilities. The distinguishing features of the model
are that the seller can control which products to recommend
depending on the arrived product and that each customer either
purchases a product or leaves the market after one transition.

We study the assortment optimization problem under this model.
Particularly, we show that this problem is generally NP-Hard even if
each product could only transit to at most two products. Despite the
complexity of the problem, we provide polynomial time algorithms for
several special cases, such as when the transition probabilities are
homogeneous with respect to the starting point, or when each product
can only transit to one other product. We also provide a tight
performance bound for revenue-ordered assortments. In addition, we
propose a compact mixed integer program formulation that can solve
this problem of large size. Through extensive numerical experiments,
we show that the proposed algorithms can solve the problem
efficiently and the obtained assortments could significantly improve
the revenue of the seller than under the Markov chain choice model.

}


\maketitle

\section{Introduction}
\label{sec:introduction}

The fast development of information technology and rapid growth of
online sales have presented great opportunities --- as well as
challenges --- for retailers to use data to increase their bottom
lines. One of the challenges is to determine which subset of
products to make available to customers in order to maximize the
expected revenue. Such a problem is often referred to as the
assortment optimization problem. Lying in the heart of the
assortment optimization problem is to identify an appropriate model
to characterize the choice behavior of consumers when facing a
subset of products. Much recent research has focused on finding
such a choice model and solving the associated assortment
optimization problem.

One of the consumer choice models that has recently gained much
attention is the Markov chain choice model proposed by
\cite{Blanchet2016}. In the Markov chain choice model, it is
assumed that customers arrive at each of the products with a certain
exogenous probability. If the product a customer arrived at is
unavailable (not in the offered assortment), then the customer will
transit to other products (including the no-purchase option) with
certain transition probabilities. And this process stops until the
customer finds an available product to purchase or leaves the
market. Many follow-up works have appeared based on the Markov chain
choice model. We will provide a more detailed literature review in
Section \ref{sec:literature}.

Although the Markov chain choice model brings a new perspective in
modeling customer choice and has been validated by empirical data,
it has two potential drawbacks. First, in the Markov chain choice
model, it is assumed that customer's transition between products
follows an exogenous process (according to the transition
probabilities). However, in practice, after a customer's arrival, it
is often the seller's recommendation that will determine the
transition of customers. This is especially true in an online
environment, in which customer's transition usually occurs when they
click a product displayed in the recommendation section associated
with an unavailable product. Second, in the Markov chain choice
model, it is implicitly assumed that customers could transit
arbitrarily many times among unavailable products, whereas in
practice, customers often have much less patience and if a product
transited to is unavailable, a customer may just leave the market
altogether. This is especially likely in the situation when the
transited product is a recommended product --- in such a case, it is
not hard to imagine that customers may become frustrated and as a
result, abandon the entire purchase session.

In this paper, motivated to overcome the above two drawbacks in the
Markov chain choice model, we propose a new choice model which we
call the Markov Chain choice model with Single Transition, or MCST
in short. Particularly, in the MCST model, customers initially
arrive at each product with a certain probability as in the Markov
chain choice model. However, if the product a customer arrived at is
not available, under the MCST model, the seller can choose a subset
of {\it available} products to recommend to the customer, and the
customer will transit to each product in the recommended set (as
well as the no-purchase option) based on a set of transition
probabilities. Note that under the MCST model, a customer will
either purchase a product or leave the market after a single
transition. We believe the MCST model exhibits the following
advantages compared to the Markov chain choice model:
\begin{enumerate}
\item Allow more flexible control for the seller in terms of which set of
products to recommend after a customer searches into an unavailable
product, which may give opportunities to achieve better revenue;
\item By requiring that the seller can only recommend available
products, the transition of customers is more realistic, which
better reflects the limited patience of customers.
\end{enumerate}

In this paper, we study the properties of the MCST model and the
associated assortment optimization problem. In particular, we make
the following theoretical contributions:
\begin{enumerate}
\item When the transition probabilities are homogeneous with respect
to the starting points, we show that the optimal assortment under
the MCST model is equivalent to that in the Markov chain choice
model. In addition, we show that a revenue-ordered assortment must
be optimal in this case and it can be solved in linear time of the
number of products;
\item When each product can only transit to one other product, we
show that the optimal assortment optimization problem under the MCST model can be
solved in polynomial time by a dynamic program;
\item When a product could transit to more than one product, we
show that the assortment optimization problem under the MCST model
is generally NP-Hard. Nevertheless, we establish a tight performance
bound for the revenue-ordered assortment in this case;
\item We establish a compact mixed integer program (MIP) formulation for
the optimal assortment optimization problem under the MCST model. The MIP can be
solved efficiently for large size problems thus is practically
useful.
\end{enumerate}

In addition to the theoretical results, we also conduct extensive
numerical experiments for the MCST model and compare the performance
with the Markov chain choice model and a recently proposed two
product nonparametric choice model \cite{Paul2016}. In
particular, we show that under the same set of arrival and
transition probabilities, it is more often that the optimal
assortment under the MCST model achieves larger revenue than that
under the Markov chain choice model or the two product nonparametric
model, showing that the additional control of the seller over the
recommended set is valuable. We also show that the optimal
assortment solved from the MCST model is more robust under model
misspecification situations.

As discussed in the beginning, sellers nowadays strive to utilize
all possible data and control to improve their sales decisions and
ultimately increase their bottom lines. The models and solutions
proposed in this paper help the seller achieve this goal by allowing
them to exert more controls in the customer choice process and
helping them obtain better decisions. Therefore, we believe that our
paper is helpful for sellers to meet the challenges presented in
this age.

The remainder of this paper is organized as follows: In Section
\ref{sec:literature}, we review related literature to our work. In
Section \ref{sec:model}, we formally state the problem studied in
this paper. In Section \ref{sec:homogeneous}, we study the problem
where the transition weights are homogeneous with respect to the
starting point. In Section \ref{sec:heterogeneous}, we investigate
the heterogeneous case. We first study a special case, and then
consider the computational complexity and algorithms for the general
case. We present some numerical experiment results in Section
\ref{sec:numerical}. Finally, we provide some concluding remarks in
Section \ref{sec:conclusion}.

\section{Related Work}
\label{sec:literature}

In this section, we review related literature to our work. At the
high level, our work is related to the fast growing research area of
choice models and the associated revenue management problems.
Particularly, our work is closely related to the recently proposed
Markov chain choice model \cite{Blanchet2016} and the
associated assortment optimization problem. In the following, we
focus our literature review on these lines of works and refer the
readers to \cite{train2009discrete} and \cite{Ozalp} for
comprehensive reviews for general research about choice models and
revenue management respectively.

The Markov chain choice model is first introduced by
\cite{Blanchet2016}. In the Markov chain choice model, there is a
line of products and each customer arrives at each product with a
certain probability. If the product is available, then the customer
will purchase it, otherwise the customer will transit to another
product or leave with certain transition probabilities. The process
terminates until the customer purchases a certain product or chooses
to leave. In \cite{Blanchet2016}, the authors show that such a model
is a good approximation to several well-studied choice models.
Moreover, they show that the assortment optimization problem under
the Markov chain choice model can be solved in polynomial time by
relating it to an optimal stopping problem on a Markov chain. There
are several follow-up works that study the Markov chain choice
model. For example, \cite{Feldman2014RevenueMU} study the network
revenue management problem under the Markov chain choice model and
present an algorithm based on linear optimization. \cite{Desir2016}
study the capacity constrained assortment optimization problem. And
\cite{desir_robust} study the robust assortment optimization under
this model. As described in the introduction, the MCST model we
propose modifies the Markov chain choice model in two important
dimensions: 1) we add a control for the seller to decide which set
of products to recommend when customers arrive at an unavailable
product and 2) we only allow a single transition of the customers.
We propose solution methods for solving the assortment optimization problem under
the modified model and show that it has the potential to increase
the revenue for the seller compared to that in the Markov chain
choice model.

Our work is also related to the growing literature that studies
assortment optimization problem. Assortment optimization is a
central problem in revenue management in which the seller selects a
subset of products to offer in order to maximize the expected
revenue. The assortment optimization problem was first considered by
\cite{Talluri2004}, in which the authors study a single resource
revenue management problem. They show that a class of
revenue-ordered assortment is optimal when customers choose
according to a multinomial logit (MNL) choice model. There are many
subsequent works studying the assortment optimization problem under
various consumer choice models. Some examples include the nested
logit model (see, e.g., \citealt{li_orl}, \citealt{Davis2014},
\citealt{gallego_nested},  \citealt{li}); the mixture of multinomial
logit (MMNL) model (see, e.g., \citealt{desir_goyal},
\citealt{rusmevichiengtong_mmnl}), the ranking-based models (see,
e.g., \citealt{Honhon2012}, \citealt{Aouad2016},
\citealt{bertsimas_misic}), the Mallows model (see e.g.,
\citealt{desir_NIPS}, \citealt{desir_mallows}) and more recently,
more complicated choice models such as the consider-then-choose
model \cite{aouad_consider_then_choose}, the MNL model with
endogenous network effects \cite{wang_wang}, choice model when
consumer searches for product information \cite{wang_sahin}, etc. In
another direction, there has also been interest in studying
assortment optimization problems under nonparametric choice models
(see \citealt{farias}, \citealt{jagabathula_nonparametric},
\citealt{Jagabathula2016}) or in a dynamic environment (see
\citealt{Rusmevichientong2010},
\citealt{Davis:2015:AOO:2846978.2846992}, \citealt{Chen2016},
\citealt{Aouad2015}). In this work, we consider the assortment
optimization problem under a new model, the MCST model. Similar to
other works that study the assortment optimization problem, we study
the complexity of the problem under different cases and propose
efficient algorithms to solve the problem. We further compare the
solution and the achieved revenue under our model with several
related models.

In the study of assortment optimization problem, an important class
of assortment is called the revenue-ordered assortment, in which the
assortment is selected according to the revenue order of each
product. Revenue-ordered assortment has been shown to be optimal in
many assortment optimization problems or has a worst-case
performance guarantee (see, e.g., \citealt{Talluri2004},
\citealt{rusmevichiengtong_robust},
\citealt{rusmevichiengtong_mmnl}, \citealt{wang_wang}). Indeed, as
shown in \cite{Berbeglia2015}, the optimal revenue-ordered assortment achieves a
worst-case ratio of the optimal revenue under a class of ``regular
discrete choice models''. In this paper, we show that the
revenue-ordered assortment is optimal under a special case of the
MCST model, and has a worst-case performance guarantee under the
general case. Furthermore, in our numerical experiments, we show
that the revenue-ordered assortment generally achieves good
performance.

One closely related work to ours is the work by \cite{Paul2016}. In
\cite{Paul2016}, the authors consider a two step choice model and a
two product nonparametric choice model. Both of those models are
based on the Markov chain choice model and only allow customers to
make a single transition. In particular, the former one only allows
homogeneous transition probabilities with respect to the starting
point, while the latter one allows heterogeneous transition
probabilities. They study the assortment optimization problem under
these models, showing that both problems are NP-Hard, while
proposing a fully polynomial time approximation scheme (FPTAS) for
the former problem and a $2$-approximation algorithm for the latter
problem. Similar models have also been considered in \cite{Kok} in
which the authors mainly consider the estimation problem and propose
heuristic algorithms for the assortment optimization problem. In
contrast to those models, in the MCST model, we allow the seller to
choose which subset of products to recommend after a customer
reaches an unavailable product instead of assuming the customer
transitions are exogenous. As we will show in our paper, such
flexibility will lead to a different optimization model for the
assortment optimization problem and the achieved revenue could be
much higher under the MCST model than in their models. Therefore,
one can view our MCST model as a modification to the models in
\cite{Paul2016}, which could lead to significant revenue improvement
for the sellers.

\section{Model}
\label{sec:model}

In this section, we introduce the MCST model (Markov chain choice
model with single transition) and the associated assortment
optimization problem. Consider a seller who manages a collection of
$n$ products denoted by ${\mathcal N} = \{1, ..., n\}$ with $0$
being the no-purchase option. We assume each product $j$ has a fixed
revenue $r_j$ (one can also view $r_j$ as the profit margin of
product $j$ if the objective is to maximize profit). Without loss of
generality, we assume that the products are sorted in the
non-increasing order based on the revenue, i.e., $r_1 \geq r_2 \geq
\cdots \geq r_n$. In the MCST model, the seller first selects an
assortment $S\subseteq {\mathcal N}$ to offer to the customers (in
other words, the seller makes available a subset of products $S$).
Then, with probability $\lambda_j \ge 0$, a customer arrives at
product $j$. (Naturally, we assume that $\sum_{j=1}^n \lambda_j =
1$.) When a customer arrives at product $j$, if product $j$ is in
the offered assortment $S$, then the customer purchases product $j$
and the seller gains revenue $r_j$. Otherwise, the seller can
recommend a subset of products within the assortment, i.e., a subset
$R_j \subseteq S$, and the customer will either select one of the
products in $R_j$ or leave the market without purchasing any
product. In particular, there is a weight $v_{ji} \ge 0$ for any
$j\in {\mathcal N}$ and $i\in {\mathcal N}\cup \{0\}$, satisfying
$\sum_{i\in {\mathcal N}\cup\{0\}} v_{ji} = 1$ for any $j$, such
that the probability a customer will transit to product $i\in R_j$
or leave the market (transit to product $i = 0$) given that the
customer arrived at product $j\notin S$ is
\begin{eqnarray*}
\rho_{ji} = \frac{v_{ji}}{\sum_{i\in R_j}v_{ji} + v_{j0}}.
\end{eqnarray*}
The assortment optimization problem under the MCST model (we will
call it the MCST problem in the following) is to find an assortment
$S \subseteq {\mathcal N}$, and for each $j \not \in S$, a
recommended subset of products $R_j \subseteq S$, such that the
expected total revenue is maximized.

Now we derive the formal expression of the MCST problem. Let
$\mathrm{Pr}_i(S,\{R_j\}_{j\not\in S})$ be the probability that $i$
is purchased by the customer when the offered assortment is $S$ and
the recommended set for $j\not\in S$ is $R_j$. For simplicity of
notation, we denote the probability by $\mathrm{Pr}_i(S,R_j)$ when
there is no confusion. We have
\begin{equation} \mathrm{Pr}_i(S,R_j) =
\begin{cases}
\lambda_i + \sum\limits_{j:j\not\in S, i\in R_j}\lambda_j\frac{v_{ji}}{\sum_{i\in R_j}v_{ji} + v_{j0}}, & \mathrm{~if~} i\in S,\\
\sum\limits_{j:j\not\in S}\lambda_j\frac{v_{j0}}{\sum_{i\in R_j}v_{ji} + v_{j0}}, & \mathrm{~if~} i = 0,\\
0, & \mathrm{~otherwise.~}
\end{cases}
\label{eq:MCST_choice_rj}
\end{equation} Let $\mathbf{R}(S,R_j)$ be the
expected revenue when the offered assortment is $S$ and the
recommended set for $j\not\in S$ is $R_j$. Then we have
\begin{alignat*}{2}
\mathbf{R}(S,R_j) = \sum_{i\in S}r_i\mathrm{Pr}_i(S,R_j) =
\sum_{i\in S}\lambda_ir_i + \sum_{j\not \in S}\sum_{i\in
R_j}\lambda_jr_i\frac{v_{ji}}{\sum_{i\in R_j}v_{ji} + v_{j0}}
\end{alignat*}
and the assortment optimization problem under the MCST model can be
written as
\begin{alignat}{2}
\max_{S\subseteq {\mathcal N}, R_j\subseteq S,\forall j\notin S}
\mathbf{R}(S,R_j)= \max_{S\subseteq {\mathcal N}, R_j\subseteq S,
\forall j\notin S} \sum_{i\in S}\lambda_ir_i + \sum_{j\not \in
S}\sum_{i\in R_j}\lambda_jr_i\frac{v_{ji}}{\sum_{i\in R_j}v_{ji} +
v_{j0}}. \label{eq:MCST}
\end{alignat}
We remark that in the case of $v_{j0} = 0$ for some $j\in S$, if the
seller recommends none of the products (i.e., $R_j = \phi$), or only
recommends products with zero transition weights ($v_{ji} = 0$),
then the ratio in (\ref{eq:MCST}) would become $0/0$. To avoid
confusion in those cases, we stipulate that $0/0 = 0$ in the rest of
our discussion. Note that this is equivalent as setting the
transition weight to the no-purchase option to an infinitesimal
value when it is zero.

By examining \eqref{eq:MCST}, if we have chosen an assortment $S$,
then the recommended set for any $j\not\in S$ can be found by
solving $\max\limits_{R_j\subseteq S}\frac{\sum_{i\in
R_j}r_iv_{ji}}{\sum_{i\in R_j}v_{ji} + v_{j0}}$. Note that this is
an assortment optimization problem under an attraction model, for which it has
been shown that a revenue-ordered assortment is optimal
\cite{Talluri2004}. Therefore, we can efficiently find the
optimal recommended set for each $j\not\in S$ if the assortment $S$
has been determined. However, finding an optimal $S$ to offer
initially is still a difficult task. As we will show in Section
\ref{sec:heterogeneous}, the problem is NP-Hard in general, and we
will develop algorithms for solving this problem in the following
sections.

Before we end this section, we comment on the relation between the
MCST model and the Markov chain choice model proposed in
\cite{Blanchet2016} and the two step choice model and the two
product nonparametric choice model proposed in \cite{Paul2016},
which are all based on similar ``transition'' ideas. Compared to
those models, a distinctive feature of the MCST model is that it
allows the seller to choose a subset of products to recommend when a
customer arrives at an unavailable product, and the subset chosen
could depend on the product the customer arrived at. In other words,
we allow the seller to control the transition of customers to some
degree rather than assuming the transition is entirely exogenous.
This is very reasonable since in practice, especially in the growing
industry of online retailing, it is often the seller's
recommendation that will determine the transition of the customers.
Moreover, we only allow the seller to recommend products that are in
the offered assortment, that is, we do not allow transition into
unavailable products. This is also very reasonable as it does not
make much sense in practice to recommend products that are
unavailable (which will make customers frustrated or even annoyed
and thus abandon the purchase session). Given the latter
restriction, customer can only transit once in the MCST model rather
than arbitrary number of times as in the Markov chain choice
model.\footnote{It can be easily shown that if one first assumes
that customers can only make a single transition, then it must be
optimal for the seller to recommend products that are in the offered
set. Therefore, the assumption that the seller only recommends
products in the offered set is equivalent to the assumption that
customers can only make a single transition. And the MCST model can
be obtained through either assumption.} Later we will show that
these features of the MCST model will lead to potentially different
optimal assortment for the seller and also higher revenue.

\section{Homogeneous Case}
\label{sec:homogeneous}

In this section, we first consider a special case of the MCST
problem in which the transition probabilities of customers are
homogeneous with respect to the starting product. More precisely, we
consider the case in which there exists $v_i \ge 0$ for $i\in
{\mathcal N}\cup \{0\}$ satisfying $\sum_{i\in {\mathcal
N}\cup\{0\}} v_i = 1$ such that $v_{ji} = v_i$ for all $i\in
{\mathcal N}\cup \{0\}$ and $j\in {\mathcal N}$. Note that this
condition ensures that the transition probability to a product $i$
given a recommended set $R$ is independent of the previous
(unavailable) product $j$ (we will call $j$ the ``starting point''
in the following discussion). In the following, we will show that
the optimal assortment under the homogeneous case of the MCST
problem (which we will call the homogeneous MCST problem in the
following) is identical to that under the Markov chain choice model
with the same parameters, has a revenue-ordered structure, and can
be solved very efficiently (in linear time of the number of
products).

To start with, we note that in the homogeneous case, since the
transition probabilities are homogeneous with respect to the
starting point, there exists an optimal solution in which the
optimal recommended sets $R_j$s are identical for all $j\notin S$.
The following theorem further shows that in this case, there is an
optimal solution in which the seller recommends all the products in
the offered sets:

\begin{lemma}\label{lem:MCST_homo}
For the homogeneous MCST problem, given any assortment $S$, let $S'
\in \mbox{arg}\max_{R\subseteq S}\frac{\sum_{i\in R}
r_iv_i}{\sum_{i\in R} v_i + v_0}$ be an optimal recommended set
corresponding to $S$ such that there does not exist $S'' \in
\mbox{arg}\max_{R\subseteq S}\frac{\sum_{i\in R} r_iv_i}{\sum_{i\in
R} v_i + v_0}$ satisfying $S' \subsetneqq S''$. Then
$\mathbf{R}\left(S, S'\right) \leq \mathbf{R}\left(S', S'\right)$.
\end{lemma}

{\noindent\bf Remark:} The condition that there does not exist $S''
\in \mbox{arg}\max_{R\subseteq S}\frac{\sum_{i\in R}
r_iv_i}{\sum_{i\in R} v_i + v_0}$ such that $S' \subsetneqq S''$ in
Lemma \ref{lem:MCST_homo} is necessary (in some sense, this
condition can be understood as when there are several recommended
sets that achieve the same revenue, we shall select the largest
one). Otherwise, there could be a $k$ with $v_k = 0$ and $r_k$ very
large such that even though $S'$ is an optimal recommended set with
respect to $S = S'\cup \{k\}$,
$\mathbf{R}(S, S') > \mathbf{R}(S', S')$. $\hfill\Box$\\

The proofs of this and all subsequent lemmas/theorems are provided
in the Appendix. Lemma \ref{lem:MCST_homo} implies that there exists
an optimal assortment $S^*$ to the homogeneous MCST problem, such
that when we choose $S^*$ to offer, we recommend the same set $S^*$
for all $j\notin S^*$. Next we show that solving the homogeneous
MCST problem is equivalent as solving the assortment optimization
problem under the Markov chain choice model with the same
parameters.

\begin{theorem}\label{th:MCST_homo_MCST}
Consider the homogeneous MCST model and the Markov chain choice
model with the same arrival probabilities $\lambda_i$ and transition
probabilities $v_{ji} = v_i$. If $S^*$ is an optimal assortment
under the Markov chain choice model, then offering $S^*$ and
recommending $S^*$ for each $j\notin S^*$ is an optimal solution to
the homogeneous MCST problem. Furthermore, the optimal values are
the same.
\end{theorem}

Lemma \ref{lem:MCST_homo} and Theorem \ref{th:MCST_homo_MCST}
suggest that we can find the optimal assortment in the homogeneous
MCST problem by solving the assortment optimization problem under
the Markov chain choice model. It is known that the latter problem
can be solved in polynomial time by relating to an optimal stopping
problem on a Markov chain \citep{Blanchet2016} or solving a linear
optimization problem \citep{Feldman2014RevenueMU}. Thus the
homogeneous MCST problem can be solved in polynomial time. In the
following, we present an iterative algorithm for solving this
problem which runs in linear time in the number of products, which
is significantly faster than existing methods for solving the optimal assortment under the Markov chain choice model. The algorithm is
motivated by the local ratio framework for the capacity constrained
assortment optimization problem under the Markov chain choice model
in \cite{Desir2016}. Furthermore, we show by the iterative algorithm
that a revenue-ordered assortment is optimal for the homogeneous
MCST problem.

Recall that the products are sorted by non-increasing order based on
their revenues. The iterative algorithm is described as follows: If
$r_1 < 0$, then it is obvious that an empty assortment is optimal.
Otherwise, we select product $1$ into the optimal assortment. Then,
we construct a new instance of an assortment optimization problem
excluding product $1$. Specifically, we keep the arrival
probabilities $\lambda_i$s unchanged and update the revenues and
transition probabilities by setting
\begin{equation} r^{(1)}_i :=
 r_i - \frac{r_1v_{1}}{v_{1} + v_{0}}, \qquad
v^{(1)}_{i}=
\begin{cases}
v_{0} + v_{1} & \mathrm{~if~} i = 0\\
v_{i} & \mathrm{otherwise}\\
\end{cases} \quad\mbox{ for all } i \in {\mathcal N}\setminus\{1\}.\label{eq:MCST_update}
\end{equation}

Let $\mathbf{R}^{(1)}(S,R)$ be the revenue of assortment $S$ under
the updated revenue $r_i^{(1)}$s and transition probabilities
$v^{(1)}_i$s, and when the recommended set for unavailable products
is $R$. We have the following lemma:

\begin{lemma}\label{lem:MCST_add}
For any $S\subseteq {\mathcal N}\setminus \{1\}$,
$\mathbf{R}(S\cup\{1\},S\cup\{1\})$ = $\mathbf{R}(\{1\},\{1\}) +
\mathbf{R}^{(1)}(S,S)$. Furthermore, when $r_1 \ge 0$,
$\mathbf{R}(S\cup\{1\},S\cup\{1\}) \ge \mathbf{R}(S, S)$.
\end{lemma}

According to Lemma \ref{lem:MCST_add}, when $r_1\ge 0$, the optimal
assortment to the homogeneous MCST problem must include product $1$.
Moreover, the optimal assortment would consist of product $1$ and
the optimal assortment to the updated instance with parameters
$r_i^{(1)}$ and $v_i^{(1)}$. To find the optimal assortment of the
updated instance, we can apply the same idea: We first check if
$r_2^{(1)} \ge 0$. If so, we include product $2$ and further update
the instance (by computing $r_i^{(2)}$ and $v_i^{(2)}$); otherwise
if $r_2^{(1)} < 0$, then the optimal assortment to the updated
instance must be an empty set. We can repeat the procedure until
there is no positive revenue product in the updated instance and the
assortment at that moment will be the optimal assortment of the MCST
problem.

In Algorithm \ref{alg:homo_iter}, we summarize the above procedure.
The theorem below it shows that Algorithm \ref{alg:homo_iter} solves
the homogeneous MCST problem in $O(n)$ time, which is significantly
faster than solving an optimal stopping problem as in
\cite{Blanchet2016} or solving an LP as in
\cite{Feldman2014RevenueMU}.

\begin{algorithm}[htb]
\caption{Iterative Algorithm for Homogeneous MCST Problem}
\label{alg:homo_iter}
\begin{algorithmic}[1]
\STATE Initialization: Let $RV_1 := r_1v_1$, $V_1 := v_0 + v_1$.
 \FOR {$i$ from $2$ to $n-1$}
 \STATE $RV_i := RV_{i-1} + r_iv_i$
 \STATE $V_i := V_{i-1} + v_i$
\ENDFOR
\FOR {$i$ from $2$ to $n$}
 \STATE   $\displaystyle r^{(i-1)}_i := r_{i} - \frac{RV_{i-1}}{V_{i-1}}$ \label{alg:homo_cal}
 \IF {$r^{(i-1)}_i < 0$}
 \STATE {\bf return} $S^* = \{1,...,i-1\}$
 \ENDIF
\ENDFOR
\STATE {\bf return} $S^* = \{1,...,n\}$.
\end{algorithmic}
\end{algorithm}

\begin{theorem}
\label{th:MCST_homo} Algorithm \ref{alg:homo_iter} finishes in
$O(n)$ time and the returned $S^*$ is an optimal assortment for the
homogeneous MCST problem.
\end{theorem}

We observe that Algorithm \ref{alg:homo_iter} always returns a
revenue-ordered assortment when it finishes, thus we have the
following corollary.

\begin{corollary}\label{coro:MCST_homo}
There is a revenue-ordered assortment that is optimal to the
homogeneous MCST problem.
\end{corollary}

In addition to obtaining efficient algorithms for solving the
homogeneous MCST problem, by Theorems \ref{th:MCST_homo_MCST},
\ref{th:MCST_homo} and Corollary \ref{coro:MCST_homo}, we know that
for the Markov chain choice model, if the transition probabilities
are homogenous with respect to the starting point, then a
revenue-ordered assortment is optimal for the assortment
optimization problem. Furthermore, there exists an algorithm that
can solve the problem in linear time. As far as we know, this is the
first time such a result is obtained which could also contribute to
the study of the Markov chain choice model.

\section{Heterogeneous Case}
\label{sec:heterogeneous}

In this section, we consider the general MCST problem in which that
transition probabilities are not necessarily homogeneous with
respect to the starting point. We study the computational complexity
of this problem as well as propose efficient algorithms for solving
it (both in special cases as well as in general cases). We start
with the following theorem stating that the problem is hard to solve
in general.

\begin{theorem}\label{th:MCST_heter_complexity}
The MCST problem is strongly NP-Hard, even if we restrict that each
product can only transit to at most two products.
\end{theorem}

The proof of Theorem \ref{th:MCST_heter_complexity} is based on a
reduction from the independent set problem and is relegated to the
Appendix. By Theorem \ref{th:MCST_heter_complexity}, it is not
possible to obtain polynomial time algorithms for solving the MCST
problem in general unless $P = NP$. Furthermore, by the strong
NP-Hardness, it is unlikely that there is an FPTAS for the MCST
problem either \cite{GJ79}. In the following, we first propose an
efficient algorithm to solve the problem when each product can only
transit to one product (Section \ref{sec:heter_one}), then establish
a worst-case performance bound for the best revenue-ordered
assortment (Section \ref{sec:heter_gen_approx}), and finally propose
a compact mixed integer program (MIP) formulation that can solve the
problem of practical size efficiently.

\subsection{Transit to One Product}
\label{sec:heter_one}

In this section, we consider a special case of the MCST problem in
which each product can only transit to one other product (in
addition to the no-purchase option). That is, for each $j\in
{\mathcal N}$, there is at most one $i \in{\mathcal N}$ such that
$v_{ji} > 0$ (in addition to $v_{j0}$). We show that the optimal
assortment in this case can be found by solving a dynamic program in
a linear time of the number of products.

\begin{figure}[htbp]
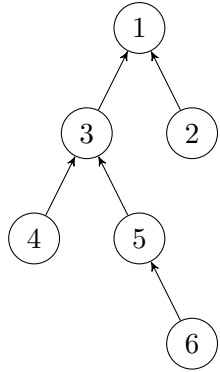

\begin{center}
\tikz[scale=.7]{ \node(v1) at (0, 0) [circle,draw] {$1$}; \node(v3)
at (-1, -2) [circle,draw] {$3$}
    edge[->,>=stealth'] (v1);
\node(v5) at (1, -2) [circle,draw] {$2$}
    edge[->,>=stealth'] (v1);
\node(v4) at (-2, -4) [circle,draw] {$4$}
    edge[->,>=stealth'] (v3);
\node(v6) at (0, -4) [circle,draw] {$5$}
    edge[->,>=stealth'] (v3);
\node(v8) at (1, -6) [circle,draw] {$6$}
    edge[->,>=stealth'] (v6);
}
\end{center}
\caption{Example of the Tree} \label{fig:tree}
\end{figure}

In the following, if $v_{ji} > 0$, then we say there is a link from
$j$ to $i$. First, if there is a link from $j$ to $k$ with $r_j \geq
r_k$, then it is easy to see that it is always optimal to include
$j$ in the assortment and the link would become useless and thus can
be safely removed. Thus, we can assume that all links from $j$ to
$k$ are such that $j
> k$, i.e., $r_j < r_k$. Therefore, we can construct a directed tree
to indicate the transition relationship between the products (an
example of such a tree is given in Figure \ref{fig:tree}). If there
is a link from $j$ to $k$, then we call $j$ a child of $k$. Without
loss of generality, we assume that product $1$ is the root of the
tree.

In the following, we show that we can solve the assortment
optimization problem in this case by solving a dynamic program. Let
$\cV(i, 1)$ denote the optimal revenue that can be obtained using
the  products in the subtree starting from product $i$ and product
$i$ is included in the assortment, and ${\cV}(i, 0)$ denote the
optimal revenue that can be obtained using the products in the
subtree starting from product $i$ and product $i$ is not included in
the assortment. We have the following recursion formula: \begin{equation}
\begin{cases}
  {\cV}(i, 1) =  \lambda_ir_i + \sum_{\{j:j \mathrm{~is~a~child~of~}i\}}\max\left\{{\cV}(j, 1), \displaystyle\frac{\lambda_jv_{ji}r_i}{v_{ji} + v_{j0}} + {\cV}(j, 0)\right\}\nonumber\\
  {\cV}(i, 0) =  \sum_{\{j:j \mathrm{~is~a~child~of~}i\}}{\cV}(j, 1)
\end{cases},
\mbox{  if $i$ is not a leaf}\end{equation} and
\begin{eqnarray*}
{\cV}(i, 1) = \lambda_ir_i,\quad {\cV}(i, 0) = 0 \quad \mbox{ if $i$
is a leaf.}
\end{eqnarray*}
To explain, ${\cV}(i, 1)$ first includes the revenue of product $i$,
which is $\lambda_ir_i$. Then for each $j$ that is a child of $i$,
there are two choices: 1) include product $j$, which leads to
${\cV}(j, 1)$ from that node and below or 2) do not include product
$j$, then there is a revenue of $\frac{\lambda_jv_{ji}r_i}{v_{ji} +
v_{j0}}$ from the transition from $j$ to $i$ and a revenue of
${\cV}(j, 0)$ from the nodes below $j$. For ${\cV}(i, 0)$, it does
not generate revenue from product $i$, and for each child $j$, the
optimal solution must include $j$ (otherwise it won't generate any
revenue because $i$ is not in the assortment). Thus ${\cV}(i, 0) =
\sum_{\{j:j \mathrm{~is~a~child~of~}i\}}{\cV}(j, 1)$.

To find the optimal assortment, we can compute the value functions
from the leaves up to the root. Eventually, ${\cV}(1,1)$ gives the
optimal value (e.g., in Figure \ref{fig:tree}, we can compute the
value of $\cV$ in the order of $6\rightarrow 5\rightarrow 4
\rightarrow 3 \rightarrow 2 \rightarrow 1$). We have the following
theorem regarding the complexity of this method.

\begin{theorem}\label{th:MCST_dp}
If each product can transit to at most one other product in addition
to the no-purchase option, then the MCST problem can be solved in
$O(n)$ time by the dynamic program described above.
\end{theorem}

\subsection{Performance Bound of Revenue-Ordered Assortments}
\label{sec:heter_gen_approx}

In Section \ref{sec:homogeneous}, we showed that a revenue-ordered
assortment is optimal to the MCST problem in the homogeneous case.
However, as shown in the beginning of this section, the MCST problem
is NP-Hard in general, thus a revenue-ordered assortment is not
necessarily optimal in general. In this section, we provide tight
performance bound for revenue-ordered assortments for the MCST
problem.

In some part, our discussion is motivated by the results in
\cite{Berbeglia2015}. In \cite{Berbeglia2015}, the authors study
the worst-case performance guarantee of the revenue-ordered
assortments under a class of choice models, which they call {\it
general discrete choice model}. Particularly, a choice model is said
to be a general discrete choice model if it satisfies the following
four conditions \cite{Berbeglia2015}:
\begin{enumerate}
\item $\mathrm{Pr}_i(S) \geq 0,\;\;\;\; \forall i\in {\mathcal N}\cup\{0\}$ and $S \subseteq {\mathcal N}$;
\item $\mathrm{Pr}_i(S) = 0,\;\;\;\; \forall i\in {\mathcal N}$ and $S \subseteq {\mathcal N}\setminus \{i\}$;
\item $\sum_{i\in S} {\mathrm {Pr}}_i(S) \leq 1,\;\;\;\; \forall S\subseteq {\mathcal N}$;
\item $\mathrm{Pr}_i(S) \leq \mathrm{Pr}_i(S'),\;\;\;\; \forall S'\subseteq S \subseteq {\mathcal N}$ and $i \in S' \cup \{0\}$.
\end{enumerate}
The first three conditions are quite basic. The last condition,
which is called the {\it regularity axiom}, says that if one shrinks
the choice set, then the probability of choosing every alternative
that remains in the choice set will increase. The following result
is proved in \cite{Berbeglia2015}:

\begin{theorem} [Theorems 3.1 and 3.2 in \citealt{Berbeglia2015}]\label{thm:berbeglia}
The best revenue-ordered assortment gives at least
$\max\left\{\frac{1}{d},
\frac{1}{1+\log{\frac{r_{\max}}{r_{\min}}}}\right\}$ fraction of the
optimal revenue under a general discrete choice model, where $d$ is
the number of distinctive revenues among the alternatives and
$r_{\max}$ and $r_{\min}$ are the highest and lowest revenues among
the alternatives.
\end{theorem}

In the MCST model, we define the choice probability of product $i$
under set $S$ as the choice probability under $S$ with the optimal
recommended sets for all $j\notin S$. That is, $\mathrm{Pr}_i(S) =
\mathrm{Pr}_i(S,R_j)$ where $R_j = \mbox{arg}\max_{R_j\subseteq S}
\frac{\sum_{i\in R_j}r_iv_{ji}}{\sum_{i\in R_j}r_iv_{ji} + v{j0}}$.
In the following, for the ease of discussion, we assume that $R_j$
is unique for every $S$. (This can be achieved by adding an
infinitesimal perturbation to the $v_{ji}$s.) The next example shows
that the MCST model in general does not satisfy the regularity
axiom, thus is not a general discrete choice model.

\begin{example}\label{example:regularity}
We consider an example with four products with revenues $r_1 = 4$,
$r_2 = 3$, $r_3 = 2$ and $r_4 = 1$. The arrival probabilities are
$\lambda_1 = \lambda_2 = \lambda_3 = \lambda_4 = 1/4$. The
transition probabilities from product $1$ to other products are
$(v_{11}, v_{12}, v_{13}, v_{14}, v_{10}) = (0 , 0 , 0 , 1/2 ,
1/2)$, and the transition probabilities from product $2$ to other
products are $(v_{21}, v_{22}, v_{23}, v_{24}, v_{20}) = (1/6 , 0 ,
1/6 , 1/3 , 1/3)$. (We do not specify the transition probabilities
for other pairs of products, as they are immaterial to the result.)
Consider the assortments $S = (1, 3, 4)$ and $S' = (3, 4)$, and the
choice probability of product $3$. It can be calculated that the
best recommended set for product $2$ under $S$ is $(1, 3)$, and the
best recommended set for product $2$ under $S'$ is $(3, 4)$. One can
calculate that the choice probability $\mathrm{Pr}_3(S) = \lambda_3
+ \lambda_2\frac{v_{23}}{v_{21} + v_{23} + v_{20}} = 0.3125$ is
strictly greater than $\mathrm{Pr}_3(S') = \lambda_3 +
\lambda_2\frac{v_{23}}{v_{23} + v_{24} + v_{20}} = 0.3$, which
violates the regularity axiom. $\hfill\Box$
\end{example}

Despite that the MCST model in general does not satisfy the
definition of general discrete choice model thus we can't directly
apply the results in \cite{Berbeglia2015}, we are still able to
prove that the optimal revenue-ordered assortment achieves a certain
fraction of the optimal revenue in the MCST problem. We have the
following theorem:

\begin{theorem}\label{th:MCST_app}
The optimal revenue-ordered assortment gives at least
$\max\left\{\frac{1}{d},\frac{1}{1+\log
\frac{r_{\max}}{r_{\min}}}\right\}$ fraction of the optimal revenue
for the MCST problem. Furthermore, the performance bound is tight in
that there exist instances such that the optimal revenue of the best
revenue-ordered assortment is $O\left(\frac{1}{d}\right)$ or
$O\left(\frac{1}{1+\log \frac{r_{\max}}{r_{\min}}}\right)$ fraction
of the optimal revenue, respectively.
\end{theorem}

\subsection{Mixed Integer Programming Formulation}
\label{sec:mip}

In this section, we show that the MCST problem can be equivalently
formulated as a mixed integer program (MIP). The MIP we formulate is
very compact with only $n$ binary variables (the rest are continuous
variables). As will be shown in Section \ref{sec:numerical}, one can
solve quite large size assortment optimization problem under the
MCST model using this MIP formulation.

To obtain the MIP formulation, we define binary variables $x_{j}$
and $y_{ji}$, $\forall i, j = 1, ..., n$, where $x_{j} = 1$
indicates that product $j$ is selected in the assortment $S$, and
$y_{ji} = 1$ indicates that product $i$ is in the recommended set
$R_j$ of product $j$. Then the optimization problem \eqref{eq:MCST}
can be formulated as the following nonlinear integer optimization
problem:
\begin{subequations}
\label{eq:MCST_mlp1}
\begin{alignat}{2}
\max\;\;\; & \sum_{j = 1}^{n}\lambda_jr_jx_j + \sum_{j=1}^n\lambda_j(1 - x_j)\frac{\sum_{i = 1}^{n}r_iv_{ji}y_{ji}}{\sum_{i = 1}^{n}v_{ji}y_{ji} + v_{j0}} \label{eq:MCST_mlp1_obj}\\
\mathrm{s.t.}\;\;\;& y_{ji} \leq x_{i},& \qquad\forall i, j = 1, \ldots , n, \label{eq:MCST_mlp1_c1}\\
& y_{ji} \leq 1 - x_{j},& \qquad\forall i, j = 1, \ldots , n, \label{eq:MCST_mlp1_c2}\\
& x_{i}, y_{ji} \in \{0, 1\}, & \qquad\forall i, j = 1, \ldots , n.
\label{eq:MCST_mlp1_c3}
\end{alignat}
\end{subequations}
In the above formulation, constraints \eqref{eq:MCST_mlp1_c1} and
\eqref{eq:MCST_mlp1_c2} guarantee that a product $i$ can be
recommended after product $j$ is visited only if product $j$ is not
included in the assortment but product $i$ is. In the following, we
assume that $v_{j0} \neq 0$ for all $j\in {\mathcal N}$ (otherwise
we can set $v_{j0}$ to be an infinitesimal number, and all the
results follow in the same way). Similar to the idea in
\cite{Davis2013}, we can rewrite \eqref{eq:MCST_mlp1} as the
following mixed integer program:
\begin{subequations}
\label{eq:MCST_mlp2}
\begin{alignat}{2}
\max\;\;\; & \sum_{j = 1}^{n}\lambda_jr_jx_j + \sum_{j=1}^n\lambda_j\sum_{i = 1}^{n}r_iz_{ji}\label{eq:MCST_mlp2_obj}\\
\mathrm{s.t.}\;\;\;& z_{ji} \leq x_{i},& \qquad\forall i, j = 1, \ldots , n, \label{eq:MCST_mlp2_c1}\\
& \sum^n_{i=1}z_{ji} + z_{j0} = 1 - x_{j},& \qquad\forall j = 1, \ldots , n, \label{eq:MCST_mlp2_c3}\\
& 0 \leq z_{ji} \leq \frac{v_{ji}}{v_{j0}}z_{j0}, & \qquad\forall i, j = 1, \ldots , n,\label{eq:MCST_mlp2_c4}\\
& x_{i}\in \{0, 1\}, & \qquad\forall i = 1, \ldots ,
n.\label{eq:MCST_mlp2_c5}
\end{alignat}
\end{subequations}

It can be easily seen that \eqref{eq:MCST_mlp2} is a relaxation for
\eqref{eq:MCST_mlp1} (by noting that one can set $z_{ji} =
(1-x_j)\frac{v_{ji}y_{ji}}{\sum_{i=1}^n v_{ji}y_{ji} + v_{j0}}$ and
$z_{j0} = (1-x_j)\frac{v_{j0}}{\sum_{i=1}^n v_{ji}y_{ji} +
v_{j0}}$). In the following, we show that the relaxation is exact.
We have the following theorem, which uses similar idea as in
\cite{Davis2013} but requires a different proof.

\begin{theorem}\label{th:mip_equiv}
If $(\bx^*, \by^*)$ is an optimal solution to \eqref{eq:MCST_mlp1},
then there exists $\bz^*$ such that $(\bx^*, \bz^*)$ is an optimal
solution to \eqref{eq:MCST_mlp2}. Conversely, if $(\bx^*, \bz^*)$ is
an optimal solution to \eqref{eq:MCST_mlp2}, then there exists
$\by^*$ such that $(\bx^*, \by^*)$ is an optimal solution to
\eqref{eq:MCST_mlp1}. Furthermore, the optimal values to
\eqref{eq:MCST_mlp1} and \eqref{eq:MCST_mlp2} are equal.
\end{theorem}

Theorem \ref{th:mip_equiv} essentially says that problem
\eqref{eq:MCST_mlp1} is equivalent to problem \eqref{eq:MCST_mlp2}.
Thus in practice, one can solve \eqref{eq:MCST_mlp2} in order to
solve the MCST problem. Note that there are only $n$ binary
variables in \eqref{eq:MCST_mlp2}. As will be shown in the next
section, the MIP formulation is very efficient for solving the MCST
problem, even for problems of practical size.

\section{Numerical Experiments}
\label{sec:numerical}

In this section, we conduct numerical experiments to test the MCST
model and the proposed algorithms. The objectives of the numerical
experiments are mainly three fold: 1) demonstrate the efficiency of
the MIP formulation for solving the MCST problem and test the
performance of the revenue-ordered assortments, 2) illustrate the
benefit of allowing the seller to choose different recommendation
sets for different products and 3) study the robustness of the
assortment decision of different models by considering several
scenarios with model misspecification. Specifically, for the second
and third parts, we will compare the performance of the optimal
assortment under the MCST model with that under the Markov chain
choice model \cite{Blanchet2016} and the two product
nonparametric choice model \cite{Paul2016}.

We start by describing the setup of the experiments. In our
experiments, we consider problem instances with
$n\in\{5,10,30,50,100,300,500\}$ products in the first part and
$n\in\{5,10,30,50\}$ products in the second and third parts (since solving large instances for the two product nonparametric choice model is challenging). Such a range should be
enough to cover most practical scenarios. We consider the following
ways to generate the parameters in the numerical experiments:
\begin{enumerate}
\item Revenue: The revenue of the products $r_j$s are generated by first generating $n$ independent
random numbers from uniform distributions on $[0, 1]$ (exponential
distributions with mean $1$, resp.) and then sorting them in
decreasing order. We use $\mathbf{UNI}$ ($\mathbf{EXP}$, resp.) to
denote the corresponding way of generating the revenues.
\item Arrival and transition probabilities: For the arrival probabilities, we generate $n$ independent random
numbers $\Lambda_j$, each from a uniform distribution on $[0, 1]$,
and set the arrival probabilities to be $\lambda_j = \Lambda_j /
\sum_{j=1}^n \Lambda_j$ for all $j$. For the transition
probabilities, we consider two methods of generating them. In the
first method, for each $j$, we generate $n+1$ independent random
numbers $\Lambda_{ji}$, each from a uniform distribution on $[0,
1]$. Then we set $v_{ji} = \Lambda_{ji}/\sum_{i =
0}^{n}\Lambda_{ji}$. We use $\mathbf{DEN}$ (meaning {\it dense}) to
denote this method. In the second method, we make the transition
matrix sparse to reflect that customer interested in one product
typically is only interested in a few other products. Particularly,
for each $\Lambda_{ji}$, with probability $1 - 1/\sqrt{n}$, we set
it to be zero and with probability $1/\sqrt{n}$, we set it to be a
uniformly distributed random number on $[0,1]$ (for
$\Lambda_{j0}$, we always set it to be a uniformly distributed
random number on $[0,1]$, meaning that the customer can always
transit to the no-purchase option). Then we set $v_{ji} =
\Lambda_{ji}/\sum_{i = 0}^{n}\Lambda_{ji}$. We use $\mathbf{SPA}$
(meaning {\it sparse}) to denote this method.
\end{enumerate}

In the presentation of the results, we will use the notation $(n,
\mathbf{A}, \mathbf{B})$ to denote the combination in which there
are $n$ products, method $\mathbf{A}$ (one of $\mathbf{UNI}$ or
$\mathbf{EXP}$) is used to generate the revenues and method
$\mathbf{B}$ (one of $\mathbf{DEN}$ or $\mathbf{SPA}$) is used to
generate the transition probabilities. All our numerical experiments
are conducted on a PC with $4$GB Memory and Intel(R) Core(TM) i5-4200U CPU (@1.6GHz) using MATLAB R2016a. And all optimization problems are solved using Gurobi 7.01
(Gurobi Optimization, Inc. 2016).

\subsection{Performance of the MIP and Revenue-Ordered Assortment}
\label{sec:numerical_MIP}

In this subsection, we first test the performance of the MIP
formulation proposed in Section \ref{sec:mip} as well as the
performance of revenue-ordered assortments. The results are shown in
Table \ref{tab:numapp}.

In Table \ref{tab:numapp}, each number is calculated by averaging
$1,000$ test instances in the corresponding combination of
parameters. The first and second columns show the average runtime
(in seconds) of the MIP formulation and the average runtime for
calculating the best revenue-ordered assortment, respectively. The
third and fourth columns show the average and worst-case performance
ratio between the best revenue-ordered assortment and the optimal
assortment, respectively. From Table \ref{tab:numapp}, we can see
that the MIP formulation is very efficient. Particularly, it can
solve problems with up to $500$ products within a matter of minutes,
and for problems that have sparse structure, it is even more
efficient. (In comparison, if one enumerates all possible
assortments, then even for problems that have $30$ products, it will
take hours to solve.) For the revenue-ordered assortment, we can see
that for the dense cases, it performs quite well, with average
performance ratio being more than $99\%$, and the performance is
especially good when the problem size is large. However, for sparse
problems, the revenue-ordered assortment performs less well,
sometimes giving rise to large performance gaps. Therefore, given
its efficiency, we suggest using the MIP formulation for practical
problems.

\begin{table}[htbp]
\begin{center}
\begin{tabular}{p{10em}<{\centering}p{7em}<{\centering}p{7em}<{\centering}p{7em}<{\centering}p{7em}<{\centering}}
                 \toprule
                  Parameters  & MIP Time (s) & RO Time (s) & RO Ave. Ratio & RO Min. Ratio \\
                  \hline
                  $(5, \mathbf{UNI}, \mathbf{DEN})$   & 0.0032  &  0.0003 &  99.72\% & 86.57\%\\
                  $(10, \mathbf{UNI}, \mathbf{DEN})$   & 0.0064  & 0.0004  & 99.76\% &  97.64\%\\
                  $(30, \mathbf{UNI}, \mathbf{DEN})$   & 0.0654  & 0.0005  & 99.93\% &  99.25\%\\
                  $(50, \mathbf{UNI}, \mathbf{DEN})$   & 0.2184  & 0.0007  & 99.96\% &  99.56\%\\
                  $(100, \mathbf{UNI}, \mathbf{DEN})$   & 1.1465  & 0.0012  & 99.98\% &  99.75\%\\
                  $(300, \mathbf{UNI}, \mathbf{DEN})$   & 20.1697  & 0.0036  & 99.99\% &  99.96\%\\
                  $(500, \mathbf{UNI}, \mathbf{DEN})$   & 70.0654  & 0.0082  & 99.99\% &  99.98\%\\
                  \hline\hline
                  $(5, \mathbf{EXP}, \mathbf{DEN})$   &    0.0032 &   0.0002 &   99.75\% &     89.36\%\\
                  $(10, \mathbf{EXP}, \mathbf{DEN})$   &    0.0062 &   0.0003 &   99.71\% &     93.58\%\\
                  $(30, \mathbf{EXP}, \mathbf{DEN})$   &    0.0538 &  0.0004  &   99.86\% &     97.88\%\\
                  $(50, \mathbf{EXP}, \mathbf{DEN})$   &    0.1580 &   0.0005 &   99.91\% &     98.15\%\\
                  $(100, \mathbf{EXP}, \mathbf{DEN})$   &    0.8802 &   0.0008 &   99.95\% &     99.07\%\\
                  $(300, \mathbf{EXP}, \mathbf{DEN})$   &   11.3842 &   0.0021 &   99.98\% &     99.79\%\\
                  $(500, \mathbf{EXP}, \mathbf{DEN})$   &   39.3100 &   0.0039 &   99.99\% &     99.92\%\\                 \hline\hline
                  $(5, \mathbf{UNI}, \mathbf{SPA})$  & 0.0034  & 0.0002  & 99.01\%  & 67.52\%\\
                  $(10, \mathbf{UNI}, \mathbf{SPA})$  & 0.0046  & 0.0004  & 98.44\%  & 88.35\%\\
                  $(30, \mathbf{UNI}, \mathbf{SPA})$  & 0.0144  & 0.0009  & 98.15\%  & 91.00\%\\
                  $(50, \mathbf{UNI}, \mathbf{SPA})$  & 0.0380  & 0.0014  & 98.40\%  & 94.09\%\\
                  $(100, \mathbf{UNI}, \mathbf{SPA})$  & 0.1570  & 0.0031  & 98.51\%  & 96.86\%\\
                  $(300, \mathbf{UNI}, \mathbf{SPA})$  & 1.8557  & 0.0153  & 98.63\%  & 97.43\%\\
                  $(500, \mathbf{UNI}, \mathbf{SPA})$  & 6.8322  & 0.0418  & 98.68\%  & 97.73\%\\
                  \hline\hline
                  $(5, \mathbf{EXP}, \mathbf{SPA})$    & 0.0034  & 0.0002  & 98.59\%  & 72.72\%\\
                  $(10, \mathbf{EXP}, \mathbf{SPA})$    & 0.0047  & 0.0004  & 97.74\%  & 86.39\%\\
                  $(30, \mathbf{EXP}, \mathbf{SPA})$    & 0.0151  & 0.0008  & 96.59\%  & 86.73\%\\
                  $(50, \mathbf{EXP}, \mathbf{SPA})$    & 0.0412  & 0.0013  & 96.09\%  & 89.38\%\\
                  $(100, \mathbf{EXP}, \mathbf{SPA})$    & 0.1642  & 0.0027  & 96.65\%  & 93.09\%\\
                  $(300, \mathbf{EXP}, \mathbf{SPA})$    & 1.9177  & 0.0111  & 97.46\%  & 94.73\%\\
                  $(500, \mathbf{EXP}, \mathbf{SPA})$    & 6.6410  & 0.0243  & 97.86\%  & 96.43\%\\
                  \bottomrule
\end{tabular}
\caption{Running Time and Performance of the MIP Formulation and Revenue-Ordered Assortments}\label{tab:numapp}
\end{center}
\end{table}

\subsection{The Benefit of Allowing Recommending Different Subsets}
\label{sec:numerical_recommendation}

In this section, we study the benefit of allowing the seller to
recommend products based on the previously visited (unavailable)
product. We will consider three models, the MCST model, the Markov
chain choice model \cite{Blanchet2016} and the two product
nonparametric model \cite{Paul2016}. Note that these three
models can be captured by the same set of parameters --- the revenue
of each product, the arrival probabilities, and the transition
probabilities. In our numerical tests, we will generate instances
according to the methods described in the beginning of Section
\ref{sec:numerical}. For each instance, we will solve the optimal
assortment under each model respectively, and calculate the optimal
revenue under each model. Note that for the MCST model, we will use
the MIP formulation in Section \ref{sec:mip} to solve for the
optimal assortment; for the Markov chain choice model, we will use
the LP formulation as described in \cite{Feldman2014RevenueMU}; and
for the two product nonparametric model, we will solve the following
binary quadratic program\footnote{In \cite{Paul2016}, the authors proposed efficient methods to solve the problem approximately. However, since we are comparing the optimal solution, we use the exact formulation. In our numerical tests, the binary quadratic program can be solved efficiently for $n$ up to 50. However, for $n$ greater than $100$, it would take several days to compute $1,000$ instances. So we only report the results up to $n = 50$.}:
\begin{equation*}\label{eq:choosy_heter_obj}
\max\limits_{\bm{x}\in\{0,
1\}^n} \;\;\; \sum_{j = 1}^{n}\lambda_jr_jx_j + \sum_{j =
1}^{n}\sum_{i = 1}^{n}\lambda_ir_jv_{ij}(1 - x_i)x_j.
\end{equation*}

The results are shown in Table \ref{tab:numl}. In Table
\ref{tab:numl}, the first section shows the comparison between the
MCST model and the Markov chain choice model. (In Table 2, we use the word {\it Choosy} to denote the two product nonparametric choice model as the authors of \cite{Paul2016} call this type of customers ``choosy customers''.) Particularly, the
first subcolumn shows the percentage of times that the optimal value
of the MCST model is larger than that of the Markov chain choice
model among the $1,000$ test instances. The second, third and
fourth subcolumns show the average, minimum and maximum ratios
between the optimal value of the MCST model and the optimal value of
the Markov chain choice model among the $1,000$ test instances,
respectively. Similarly, the second section shows the comparison
between the MCST model and the two product nonparametric choice
model. Particularly, the first subcolumn shows the percentage of
times that the optimal value of the MCST model is larger than that
of the two product nonparametric choice model among the $1,000$ test
instances. The second, third and fourth subcolumns show the average,
the minimum and the maximum ratios between the optimal value of the
MCST model and the optimal value of the two product nonparametric
choice model among the $1,000$ test instances, respectively.

\begin{table}[htbp]
\begin{center}
\begin{tabular}{p{8em}<{\centering}p{5em}<{\centering}p{3em}<{\centering}p{3em}<{\centering}p{3em}<{\centering}p{5em}<{\centering}p{3em}<{\centering}p{3em}<{\centering}p{3em}<{\centering}}
                  \toprule
                    &  \multicolumn{4}{c}{MCST/Markov} & \multicolumn{4}{c}{MCST/Choosy}\\
                    \cmidrule(l){2-5}\cmidrule(l){6-9}
                   {Parameters}  & {MCST$\geq$ Markov} & Ave. Ratio & Min. Ratio & Max. Ratio & {MCST$\geq$ Choosy} & Ave. Ratio & Min. Ratio & Max. Ratio\\
                  \hline
                  {\small$(5, \mathbf{UNI}, \mathbf{DEN})$} &70.00\%  & 1.0123  & 0.7758  & 1.2405  & 100\%  & 1.1582  & 1.0000 & 2.5825\\
                  {\small$(10, \mathbf{UNI}, \mathbf{DEN})$} & 75.80\%  & 1.0142  & 0.9384  & 1.2210  & 100\%  & 1.2350  & 1.0149 &  1.9259\\
                  {\small$(30, \mathbf{UNI}, \mathbf{DEN})$} & 85.80\%  & 1.0098  & 0.9811  & 1.0494  & 100\%  & 1.3674  & 1.1456 &  1.8206\\
                  {\small$(50, \mathbf{UNI}, \mathbf{DEN})$} & 93.20\%  & 1.0089  & 0.9906  & 1.0334  & 100\%  & 1.4268  & 1.2276  & 1.8000\\
                   \hline\hline
                  {\small$(5, \mathbf{EXP}, \mathbf{DEN})$} & 61.20\%  &  1.0152  &  0.7194  &  1.4270  &  100\%  &  1.3083  &  1.0000 &  3.0025\\
                  {\small$(10, \mathbf{EXP}, \mathbf{DEN})$} & 69.90\%  &  1.0241  &  0.8158  &  1.2339  &  100\%  &  1.5497  &  1.0622  &  3.5257\\
                  {\small$(30, \mathbf{EXP}, \mathbf{DEN})$} & 87.60\%  &  1.0316  &  0.8411  &  1.1890  &  100\%  &  2.0562  &  1.4099  &  3.5480\\
                  {\small$(50, \mathbf{EXP}, \mathbf{DEN})$} & 93.40\%  &  1.0341  &  0.9181  &  1.1855  &  100\%  &  2.3341  &  1.5786  &  4.4266\\
                   \hline\hline
                  {\small$(5, \mathbf{UNI}, \mathbf{SPA})$}  & 79.80\%  &  1.0072  &  0.6663  &  1.6459  &  100\%  &  1.0570  &  1.0000  &  1.9723\\
                  {\small$(10, \mathbf{UNI}, \mathbf{SPA})$} & 62.50\%  &  1.0058  &  0.8305  &  1.2445  &  100\%  &  1.0770  &  1.0000 &  1.4575\\
                  {\small$(30, \mathbf{UNI}, \mathbf{SPA})$}   & 66.50\%  &  1.0078  &  0.9138  &  1.1038  &  100\%  &  1.1298  &  1.0167  &  1.3458\\
                  {\small$(50, \mathbf{UNI}, \mathbf{SPA})$}  & 65.20\%  &  1.0059  &  0.9570  &  1.0607  &  100\%  &  1.1573  &  1.0532  &  1.3249\\
                  \hline\hline
                   {\small$(5, \mathbf{EXP}, \mathbf{SPA})$} &  71.30\%  &1.0097  &0.6601  &1.7145  &100\%  &1.1055  &1.0000 &2.2900\\
                   {\small$(10, \mathbf{EXP}, \mathbf{SPA})$} & 56.20\%  &1.0081  &0.6320  &1.3988  &100\%  &1.1582  &1.0000  &2.1188\\
                   {\small$(30, \mathbf{EXP}, \mathbf{SPA})$} & 59.30\%  &1.0138  &0.7800  &1.3191  &100\%  &1.3032  &1.0879  &2.0047\\
                   {\small$(50, \mathbf{EXP}, \mathbf{SPA})$} & 61.80\%  &1.0131  &0.8690  &1.2225  &100\%  &1.3809  &1.1786  &1.9767\\                  \bottomrule
\end{tabular}
\caption{Comparison between the MCST and Other
Models}\label{tab:numl}
\end{center}
\end{table}

There are two forces that drive the comparison between the MCST and
the Markov chain choice model. On the one hand, the MCST allows the
seller to recommend different subset of products depending on the
previously visited product, which would increase the optimal
revenue. On the other hand, the MCST only allows a single transition
of the customer, which would decrease the optimal revenue. As shown
in Table \ref{tab:numl}, the two forces roughly balance on average
in terms of the performance ratio, with MCST generating higher
optimal value more frequently. Here we comment that the advantage of
the MCST requires the seller to take action (of recommending
different subsets) while the disadvantage is more of a modeling
assumption. Particularly, customers who will transit multiple times
will still transit multiple times in the MCST model, which will make
the disadvantage less significant. This becomes more evident when we
compare the MCST model and the two product nonparametric model. In
this comparison, only the advantage of the MCST model exists
(customers can only transit once under both models). From Table
\ref{tab:numl}, we can see that the benefit of allowing the seller
to recommend different subsets based on the starting point is
significant, with an average of $10\%-40\%$ revenue increment in
many parameter settings and over $130\%$ in some instances.
Therefore, the results indicate that allowing the seller to
recommend different subset of products based on the arrival product
could be a revenue driver in practice.

\subsection{Robustness of the Results}
\label{sec:numerical robust}

In this section, we study the robustness of the assortment generated from
the MCST model (as well as other models) by considering several
model misspecification scenarios. Particularly, we investigate the
loss of revenue if we use the optimal assortment solved from one
model but the customer behaves according to another model. Since the
optimal assortment from the MCST model is clearly the best if the
underlying model is indeed the MCST model, we will mainly consider
the cases in which the underlying model is either the Markov chain
choice model or the two product nonparametric model.

\begin{table}[htbp]
\begin{center}
\begin{tabular}{p{8em}<{\centering}p{3.5em}<{\centering}p{3.5em}<{\centering}p{3.5em}<{\centering}p{3.5em}<{\centering}p{3.5em}<{\centering}p{3.5em}<{\centering}p{3.5em}<{\centering}p{3.5em}<{\centering}}
                  \toprule
                     &  \multicolumn{4}{c}{Average Ratio} & \multicolumn{4}{c}{Minimum Ratio} \\
                  \cmidrule(l){2-5}\cmidrule(l){6-9}
                   {Parameters}  & \multicolumn{2}{c}{True Markov} & \multicolumn{2}{c}{True Choosy} & \multicolumn{2}{c}{True Markov} & \multicolumn{2}{c}{True Choosy} \\
                    \cmidrule(l){2-3}\cmidrule(l){4-5}\cmidrule(l){6-7}\cmidrule(l){8-9}
                    &  MCST & Choosy & MCST & Markov & MCST & Choosy & MCST & Markov\\
                \hline
                  {\small$(5, \mathbf{UNI}, \mathbf{DEN})$} & 0.9940  &0.9408  &0.9229  &0.9087  &0.7314  &0.5245  &0.4579  &0.4219\\
                  {\small$(10, \mathbf{UNI}, \mathbf{DEN})$} & 0.9933  &0.8950  &0.8437  &0.8106  &0.8139  &0.6508  &0.4740  &0.3782\\
                  {\small$(30, \mathbf{UNI}, \mathbf{DEN})$} & 0.9947  &0.8160  &0.6605  &0.6048  &0.9268  &0.6392  &0.3933  &0.3068\\
                  {\small$(50, \mathbf{UNI}, \mathbf{DEN})$} & 0.9953  &0.7858  &0.5738  &0.5117  &0.9554  &0.6685  &0.3118  &0.2697\\
                \hline\hline
                  {\small$(5, \mathbf{EXP}, \mathbf{DEN})$} & 0.9875  &0.8908  &0.9129  &0.8760  &0.5986  &0.4628  &0.4896  &0.2620\\
                  {\small$(10, \mathbf{EXP}, \mathbf{DEN})$} & 0.9755  &0.7862  &0.8015  &0.7379  &0.7135  &0.3699  &0.3989  &0.2526\\
                  {\small$(30, \mathbf{EXP}, \mathbf{DEN})$} & 0.9669  &0.6170  &0.5781  &0.4744  &0.7364  &0.3601  &0.2664  &0.1434\\
                  {\small$(50, \mathbf{EXP}, \mathbf{DEN})$} & 0.9660  &0.5475  &0.4734  &0.3729  &0.7039  &0.3142  &0.2633  &0.0856\\
                \hline\hline
                {\small$(5, \mathbf{UNI}, \mathbf{SPA})$} &0.9926  &0.9799  &0.9845  &0.9733  &0.6568  &0.6384  &0.7092  &0.4048\\
                {\small$(10, \mathbf{UNI}, \mathbf{SPA})$} & 0.9860  &0.9657  &0.9670  &0.9492  &0.8242  &0.7850  &0.7073  &0.6487\\
                {\small$(30, \mathbf{UNI}, \mathbf{SPA})$} & 0.9824  &0.9350  &0.9311  &0.8953  &0.8913  &0.7735  &0.7875  &0.6723\\
                {\small$(50, \mathbf{UNI}, \mathbf{SPA})$} & 0.9830  &0.9171  &0.9083  &0.8632  &0.9184  &0.8257  &0.7884  &0.6740\\
                \hline\hline
                {\small$(5, \mathbf{EXP}, \mathbf{SPA})$} & 0.9768  &0.9606  &0.9785  &0.9531  &0.6490  &0.6500  &0.5244  &0.3588\\
                {\small$(10, \mathbf{EXP}, \mathbf{SPA})$} & 0.9593  &0.9261  &0.9673  &0.9180  &0.6489  &0.5516  &0.6971  &0.4316\\
                 {\small$(30, \mathbf{EXP}, \mathbf{SPA})$} & 0.9314  &0.8540  &0.9210  &0.8224  &0.6245  &0.5817  &0.7602  &0.4510\\
                {\small$(50, \mathbf{EXP}, \mathbf{SPA})$} & 0.9244  &0.8185  &0.8940  &0.7705  &0.7179  &0.6065  &0.7424  &0.3766\\
             \bottomrule
\end{tabular}
\caption{Robustness of Different Models}\label{tab:nummis}
\end{center}
\end{table}

Again, we will consider the experimental setups as described in the
beginning of Section \ref{sec:numerical}. The results are shown in
Table \ref{tab:nummis}. In Table \ref{tab:nummis}, we show the average and worst-case performance ratios when the underlying model is one model and one
uses the optimal assortment solved by assuming another model. For
example, the True Markov-MCST column under the Average Ratio section
shows the average performance ratio of using the optimal assortment
solved from the MCST model when the true underlying model is the
Markov chain choice model versus using the optimal assortment of
the true Markov model. From Table \ref{tab:nummis}, we can see that
our model is more robust than the other models. In particular, when the underlying model is misspecified, using the assortment solved from the MCST model usually results in better performance than using the third model. This suggests using the optimal assortment under the MCST model is a robust choice.

\section{Conclusion}\label{sec:conclusion}

In this paper, we considered a Markov chain choice model with single
transition. In particular, the seller can choose a recommendation
set based on the visited product. We proposed algorithms for solving
the assortment optimization problem under this model. We showed that
this problem can be solved in polynomial time under several special
cases and proposed an efficient MIP formulation for solving the
general problem. We demonstrated via numerical experiments that the
model has potential for increasing the revenue of the seller.

There are several future directions for research. First,
theoretically, it would be interesting to see if there are other
cases in which the assortment optimization problem under the MCST
model can be solved in polynomial time. Second, it would be
important to study the estimation problem under this model. Last, it
would be very interesting to validate this model using real data and
see how effective the optimal assortment of the MCST model is in
practice.

\bibliographystyle{ormsv080}
\bibliography{ref_assort}

\section*{Appendix}\label{app:proofs}
{\noindent\bf Proof of Lemma \ref{lem:MCST_homo}.} The lemma is
trivial if $S = S'$. In the following, we consider the case when
$S\setminus S' \neq \phi$. Let $k\in S\setminus S'$, we show that
$\mathbf{R}\left(S, S'\right) \leq \mathbf{R}\left(S\setminus\{k\},
S'\right)$. Since $S'$ is an optimal recommended set under $S$, we
have
\begin{equation}\label{eq:lemma1proof}
\frac{\sum_{i\in S'}r_iv_{i}}{\sum_{i\in S'}v_{i} +
v_{0}} \geq \frac{\sum_{i\in S'\cup\{k\}}r_iv_{i}}{\sum_{i\in
S'\cup\{k\}}v_{i} + v_{0}}.
\end{equation}
By the assumption in the lemma, we must have $v_k \neq 0$, otherwise
$S'' = S' \cup \{k\} \in \mbox{arg}\max_{R\subseteq
S}\frac{\sum_{i\in R} r_iv_i}{\sum_{i\in R} v_i + v_0}$ which
contradicts with the assumption. Combine (\ref{eq:lemma1proof}) and
that $v_k \neq 0$, we have $\displaystyle r_{k} \leq
\frac{\sum_{i\in S'}r_iv_{i}}{\sum_{i\in S'}v_{i} + v_{0}}$. Thus
$\mathbf{R}\left(S\setminus\{k\}, S'\right) - \mathbf{R}\left(S,
S'\right) = \lambda_k \frac{\sum_{i\in S'}r_iv_{i}}{\sum_{i\in
S'}v_{i} + v_{0}} - \lambda_k r_{k} \ge 0$. Repeating the above
procedure, we obtain $\mathbf{R}\left(S, S'\right) \leq
\mathbf{R}\left(S', S'\right)$ and the lemma holds. $\hfill\Box$\\

{\noindent\bf Proof of Theorem \ref{th:MCST_homo_MCST}.} According
to Lemma \ref{lem:MCST_homo}, there exists an optimal solution to
the homogeneous MCST problem in which the recommended set is the
entire set of available products. Therefore, we only need to
consider the best $\mathbf{R}(S, S)$ in order to solve the
homogeneous MCST problem. In the following, we show that for any
assortment $S$, $\mathbf{R}\left(S, S\right)$ is equal to the
revenue of assortment $S$ under the Markov chain choice model. To
show that, it suffices to show that the choice probabilities of each
product are identical for the two models. By
\eqref{eq:MCST_choice_rj}, the choice probability for each $i \in S$
in the homogeneous MCST problem is $\mathrm{Pr}_i(S,S) = \lambda_i +
\sum\limits_{j:j\not\in S}\lambda_j\frac{v_{i}}{\sum_{i\in S}v_{i} +
v_{0}}$. Now we consider the choice probability under the Markov
chain choice model. From Theorem 2.1 in \cite{Blanchet2016}, the
choice probability $\mathrm{Pr}^{M}_i(S) = \lambda_i +
(\bm{\lambda}(\bar{S}))^{T}({\mathcal I} - \rho(\bar{S},
\bar{S}))^{-1}\rho(\bar{S}, \{i\})$, where $\bar{S} = {\mathcal N}
\setminus S$, $\bm{\lambda}(\bar{S})$ is the vector of arrival
probabilities in $\bar{S}$, ${\mathcal I}$ is the identity matrix
and $\rho(S, S')$ is the transition matrix from products in $S$ to
products in $S'$. Note that $\rho(\bar{S}, \{i\}) = (v_{i}, ...,
v_{i})^{T}$ and $\sum_{i\in {\mathcal N} \cup \{0\}}v_i = 1$, it can
be verified that $({\mathcal I} - \rho(\bar{S},
\bar{S}))^{-1}\rho(\bar{S}, \{i\}) = \left({v_i}/(v_{0} + \sum_{i\in
S}v_i), ..., {v_i}/(v_{0} + \sum_{i\in S}v_i)\right)^{T}$, and hence
$\mathrm{Pr}^{M}_i(S) = \lambda_i + \sum\limits_{j:j\not\in
S}\lambda_j\frac{v_{i}}{\sum_{i\in S}v_{i} + v_{0}} =
\mathrm{Pr}_i(S,S)$. Thus the theorem holds. $\hfill\Box$\\

{\noindent\bf Proof of Lemma \ref{lem:MCST_add}.} We first prove the first part. By simple algebra, we have
\begin{alignat}{2}
& \hspace{-1cm} \mathbf{R}(S\cup\{1\},S\cup\{1\})  \nonumber\\
 = & \sum_{j\in S\cup
\{1\}}\lambda_jr_j + \sum_{j\in {\mathcal N}\setminus (S\cup
\{1\})}\lambda_j\frac{\sum_{i\in
S\cup\{1\}}r_iv_{i}}{\sum_{i\in S\cup\{1\}}v_{i} + v_{0}}\nonumber\\
=&\sum_{j\in S\cup \{1\}}\lambda_jr_j + \sum_{j\in {\mathcal
N}\setminus (S\cup \{1\})}\lambda_j\frac{\sum_{i\in S}r_iv_{i} +
r_1v_1}{\sum_{i\in S}v_{i} + v_{1} + v_{0}} \nonumber\\
 =& \sum_{j\in S\cup \{1\}}\lambda_jr_j + \sum_{j\in {\mathcal
N}\setminus (S\cup \{1\})}\lambda_j\left(\frac{r_{1}v_{1}}{v_{1} +
v_{0}} +\frac{\sum_{i\in S}\left(r_i - \frac{r_1v_{1}}{v_{1} +
v_{0}}\right)v_{i}}{\sum_{i\in S}v_{i} + v_{1} + v_{0}} \right)
\nonumber\\
=& \lambda_1r_1 + \sum_{j\in S}\lambda_j\frac{r_{1}v_{1}}{v_{1} +
v_{0}} + \sum_{j\in {\mathcal N}\setminus (S\cup
\{1\})}\lambda_j\frac{r_{1}v_{1}}{v_{1} + v_{0}}\nonumber\\
& \quad\quad+ \sum_{j\in S}\lambda_j\left(r_j -
\frac{r_1v_{1}}{v_{1} + v_{0}}\right) + \sum_{j\in N\setminus (S\cup
\{1\})}\lambda_j\frac{\sum_{i\in S}\left(r_i - \frac{r_1v_{1}}{v_{1}
+ v_{0}}\right)v_{i}}{\sum_{i\in S}v_{i} + v_{1} + v_{0}} \nonumber\\
=&  \lambda_1r_1 + \sum_{j\in {\mathcal
N}\setminus\{1\}}\lambda_j\frac{r_{1}v_{1}}{v_{1} + v_{0}} +
\sum_{j\in S}\lambda_jr^{(1)}_j + \sum_{j\in ({\mathcal
N}\setminus\{1\})\setminus S}\lambda_j\frac{\sum_{i\in
S}r^{(1)}_iv^{(1)}_{i}}{\sum_{i\in
S}v^{(1)}_{i} + v^{(1)}_{0}} \nonumber\\
= & \mathbf{R}(\{1\},\{1\}) + \mathbf{R}^{(1)}(S,S).\nonumber
\end{alignat}
Next, we prove the second part. We have when $r_1\ge 0$
\begin{alignat}{2}
\mathbf{R}\left(S\cup \{1\},S\cup \{1\}\right) & = \sum_{j\in S \cup \{1\}}\lambda_jr_j + \sum_{j\not \in S \cup \{1\}}\lambda_j \frac{\sum_{i\in S\cup\{1\}}r_iv_{i}}{\sum_{i\in S\cup\{1\}}v_{i} + v_{0}} \nonumber\\
& = \sum_{j\in S}\lambda_jr_j + \sum_{j\not \in S}\lambda_j \frac{\sum_{i\in S\cup\{1\}}r_iv_{i}}{\sum_{i\in S\cup\{1\}}v_{i} + v_{0}} + \lambda_1r_1 - \lambda_1\frac{\sum_{i\in S\cup\{1\}}r_iv_{i}}{\sum_{i\in S\cup\{1\}}v_{i} + v_{0}} \nonumber\\
&\geq \sum_{j\in S}\lambda_jr_j + \sum_{j\not \in S}\lambda_j \frac{\sum_{i\in S}r_iv_{i}}{\sum_{i\in S}v_{i} + v_{0}} + \lambda_1r_1 - \lambda_1\frac{\sum_{i\in S\cup\{1\}}r_iv_{i}}{\sum_{i\in S\cup\{1\}}v_{i} + v_{0}}  \nonumber\\
&= \mathbf{R}(S,S) + \lambda_1\left(r_1 - \frac{\sum_{i\in
S\cup\{1\}}r_iv_{i}}{\sum_{i\in S\cup\{1\}}v_{i} + v_{0}}
\right)\geq \mathbf{R}(S,S), \nonumber
\end{alignat}
where the first inequality holds because when $r_1 \ge 0$,
$\displaystyle\frac{\sum_{i\in S\cup\{1\}}r_iv_{i}}{\sum_{i\in
S\cup\{1\}}v_{i} + v_{0}} \geq \frac{\sum_{i\in
S}r_iv_{i}}{\sum_{i\in S}v_{i} + v_{0}}$  and the last inequality
holds because when $r_1 \ge 0$, $r_1 \ge \frac{\sum_{i\in
S\cup\{1\}}r_iv_{i}}{\sum_{i\in S\cup\{1\}}v_{i} + v_{0}}$. Thus the
lemma is proved. $\hfill\Box$\\

{\noindent\bf Proof of Theorem \ref{th:MCST_homo}.} It is clear that
the computational complexity of Algorithm \ref{alg:homo_iter} is
$O(n)$. It remains to show that Step \ref{alg:homo_cal} actually
computes the value of updated revenue in \eqref{eq:MCST_update} in
each iteration. We show this by induction. In fact, we will show a
stronger statement that $r_j^{(i-1)} = r_j -
\frac{RV_{i-1}}{V_{i-1}}$ gives the updated revenue for all $j\ge i$
in $i$th iteration. (Note that products $1, ..., i-1$ have been
selected before $i$th iteration and there is no need to update
them.) First, the result holds for $i = 2$ by definition in
(\ref{eq:MCST_update}). Now suppose the induction assumption holds
in the $k$th iteration, i.e., $r^{(k-1)}_{i} = r_{i} -
\frac{RV_{k-1}}{V_{k-1}}$ for all $i \geq k$. Then, for the
$(k+1)$th iteration, by \eqref{eq:MCST_update}, the updated revenue
of product $i \geq k+1$ is
\begin{eqnarray*}
r^{(k)}_{i} & = & r^{(k-1)}_{i} -
\frac{r^{(k-1)}_{k}v^{(k-1)}_{k}}{v^{(k-1)}_0 + v_k}\\
& = & r_i - \frac{RV_{k-1}}{V_{k-1}} - \frac{\left(r_{k} -
\frac{RV_{k-1}}{V_{k-1}}\right)v_k}{\sum_{j = 1}^{k-1}v_j + v_{k}}\\
& = & r_i - r_k + \left(r_k - \frac{RV_{k-1}}{V_{k-1}}\right)\cdot \frac{V_{k-1}}{V_k}\\
& = & r_i - \frac{r_kv_k}{V_k} -
\frac{RV_{k-1}}{V_k}\\
& = & r_i - \frac{RV_{k}}{V_k}
\end{eqnarray*}
where $v_0^{(k-1)} = \sum_{j=1}^{k-1} v_j$ and $v_k^{(k-1)} = v_k$.
Therefore, the induction assumption holds for $i = k+1$. Thus, by
the induction principle and Lemma \ref{lem:MCST_add}, Theorem
\ref{th:MCST_homo} holds. $\hfill\Box$

{\noindent\bf Proof of Theorem \ref{th:MCST_heter_complexity}.} We
prove the theorem by a reduction from the unweighted independent set
problem, which is known to be NP-Hard in the strong sense
\citep{GJ79}: \vspace{.3cm}

\noindent{\sc Independent Set Problem:} Given a graph $G = (V, E)$,
and an integer $0 < k \leq n$, the problem is to decide if $G$
contains a subset $I\subseteq V$ such that $|I| \geq k$ and no two
vertices in $I$ are joined by an edge in $E$. \vspace{.3cm}

\begin{figure}[htbp]
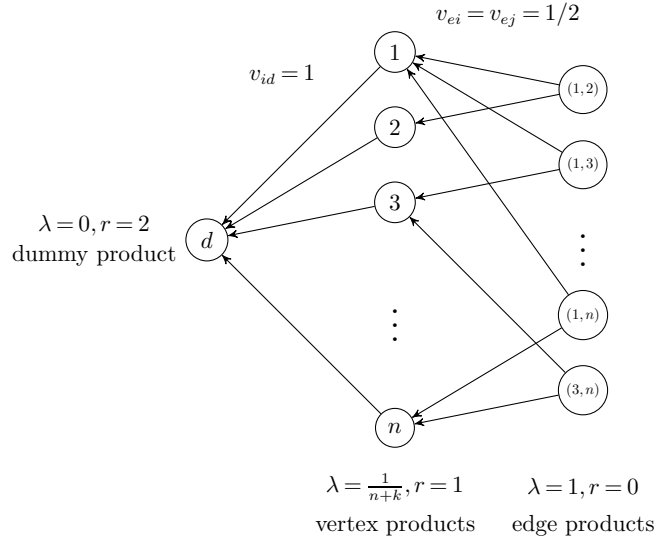

\begin{center}
\tikz{ \node(d) at (0, 3) [circle,scale=.8,draw] {$d$}; \node at
(-1.5, 3.2) [scale=.8] {$\lambda = 0,r = 2$}; \node at (-1.5, 2.8)
[scale=.8] {dummy product}; \node(v1) at (2.5, 5.5)
[circle,scale=.8,draw] {${1}$}
    edge[->,>=stealth'] (d);
\node(v2) at (2.5, 4.5) [circle,scale=.8,draw] {${2}$}
    edge[->,>=stealth'] (d);
\node(v3) at (2.5, 3.5) [circle,scale=.8,draw] {${3}$}
    edge[->,>=stealth'] (d);
\node at (2.5,2) [scale=1.2] {$\vdots$}; \node(vn) at (2.5, 0.5)
[circle,scale=.8,draw] {${n}$}
    edge[->,>=stealth'] (d);
\node(e12) at (5, 5) [circle,scale=.5,draw] {${(1,2)}$}
    edge[->,>=stealth'] (v1)
    edge[->,>=stealth'] (v2);
\node(e13) at (5, 4) [circle,scale=.5,draw] {${(1,3)}$}
    edge[->,>=stealth'] (v1)
    edge[->,>=stealth'] (v3);
\node at (5,3) [scale=1.2] {$\vdots$}; \node(e1n) at (5, 2)
[circle,scale=.5,draw] {${(1,n)}$}
    edge[->,>=stealth'] (v1)
    edge[->,>=stealth'] (vn);
\node(e3n) at (5, 1) [circle,scale=.5,draw] {${(3,n)}$}
    edge[->,>=stealth'] (v3)
    edge[->,>=stealth'] (vn);
\node at (4, 6) [scale=.8] {$v_{ei} = v_{ej} = 1/2$}; \node at (5,
-0.3) [scale=.8] {$\lambda = 1,r = 0$}; \node at (5, -0.8)
[scale=.8] {edge products}; \node at (2.5, -0.3) [scale=.8]
{$\lambda = \frac{1}{n+k},r = 1$}; \node at (2.5, -0.8) [scale=.8]
{vertex products}; \node at (1, 5.2) [scale=.8] {$v_{id} = 1$}; }
\end{center}
\caption{Illustration of the Reduction from the Independent Set
Problem} \label{fig:MCST_hard}
\end{figure}

Given an instance of the independent set problem with $G = (V,E)$,
we construct an instance of the MCST problem with $m + n + 1$
products, where $n = |V|$ and $m = |E|$. For ease of exposition, we
assume that the arrival probabilities are unnormalized, i.e., they
might be greater than 1. This is without loss of generality since we
can normalize them by dividing their sum. The products are divided
into three classes. For each edge $e \in E$, we create a product,
called the edge product, each has arrival probability $1$ and
revenue $0$. For each vertex $i\in V$, we create a product, called
the vertex product, each has arrival probability $\frac{1}{n+k}$ and
revenue $1$. There is also a dummy product $d$, which has arrival
probability $0$ and revenue $2$. The transition weights between the
products are defined as follows. The edge products can transit to
the vertex products. Specifically, if $e = (i, j)$, then the
transition weights are $v_{ei} = v_{ej} = 1/2$, and the weights are
$0$ otherwise (including the no-purchase option). Each vertex
product can only transit to the dummy product, and the weights are
$v_{id} = 1$ and $v_{i0} = 0$ for any $i\in V$. It can be seen that
in the construction, the largest value is polynomially bounded by
the size of the input instance (thus it is a polynomial-time reduction). We illustrate the construction in
Figure \ref{fig:MCST_hard}.

In the following, we show that there exists an independent set of
size at least $k$ if and only if the optimal assortment to the MCST
problem has a revenue of at least $m + 1$. If this holds, then by
the strong NP-Hardness of the independent set problem, we know that
the MCST problem is also strongly NP-Hard (and note that in the
construction each product only transits to at most two products,
thus the theorem holds).

First, we prove that if there exists an independent set of size $k$,
then there exists an assortment in the MCST problem that gives a
revenue of at least $m + 1$. By the construction of the instance, it
is obvious that an optimal assortment must include the dummy product
which has zero arrival probability, but not any edge product which
has zero revenue. Let $U$ be the subset of vertex products that are
selected into the assortment, and $\overline{U}$ be those that are
not in the assortment. Denote ${\bf R}_E$ as the revenue of the edge
products (obtained by transiting to the products in $U$), ${\bf
R}_U$ as the revenue of products $U$ and ${\bf R}_{\overline{U}}$ as
the revenue of products $\overline{U}$ (obtained by transiting to
the dummy product). Note that the revenue of the dummy product
itself is 0 since its arrival probability is 0. Thus, the total
revenue can be represented as ${\bf R} = {\bf R}_E + {\bf R}_U +
{\bf R}_{\overline{U}}$.

Suppose that there exists an independent set $I$ of size at most
$k$. We consider the assortment by selecting the vertex product into
$U$ if and only if its corresponding vertex is not in the
independent set $I$. For each edge product, the revenue is $1$ if at
least one of its endpoints is in $U$ (regardless of which one or
both), and the revenue is $0$ if none of its endpoints is in $U$
(recall we assumed in the remark after (\ref{eq:MCST}) that the
revenue is $0$ in such cases, and such assumption is without loss of
generality because we can choose the transition probability to the
no-purchase option as an infinitesimal value instead of $0$ in the
construction). Since $I$ is an independent set, at least one of the
endpoint is not in $I$ for each edge, thus it follows that ${\bf
R}_E = m$. Furthermore, the revenue of all the selected vertex
products is ${\bf R}_U = \frac{|U|}{n+k} = \frac{n - |I|}{n+k}$, and
the remaining revenue is ${\bf R}_{\overline{U}} =
\frac{2|I|}{n+k}$. Therefore, there exists an assortment which has
revenue ${\bf R} = {\bf R}_E + {\bf R}_U + {\bf R}_{\overline{U}} =
m + \frac{n+|I|}{n+k} \geq m+1$ since $|I|\geq k$.

Now we prove the other direction. Suppose that there exists an
assortment which has revenue at least ${\bf R} \geq m + 1$. We
select the vertices into the independent set $I$ if and only if the
corresponding vertex product does not belong to the assortment $U$.
We first claim that ${\bf R}_E \geq m$, which implies ${\bf R}_E =
m$ since each edge product can contribute to the revenue by at most
1. Assume that ${\bf R}_E \leq m - 1$. By the arguments above, we
have ${\bf R}_U + {\bf R}_{\overline{U}} = \frac{n + |I|}{n+k}$.
Note that $k > 0$ and $|I|\leq n$, it follows that ${\bf R} = {\bf
R}_E + {\bf R}_U + {\bf R}_{\overline{U}} \leq m - 1 + \frac{n +
|I|}{n+k} < m - 1 + \frac{2n}{n} = m + 1$, which contradicts with
the assumption that ${\bf R} \geq m + 1$. Therefore ${\bf R}_E = m$
and each edge has at least one endpoint in $U$. In other words, the
vertex set $I$ constitutes an independent set. Furthermore, we have
${\bf R}_U + {\bf R}_{\overline{U}} = \frac{n + |I|}{n+k} \geq 1$,
and thus $|I| \geq k$. Therefore, the instance has an independent
set with size at least $k$. Thus, the theorem holds. $\hfill\Box$\\

{\noindent\bf Proof of Theorem \ref{th:MCST_dp}.}  The correctness
of the DP directly follows from the recursion formula. Now we briefly
discuss the computational complexity. We note that in the dynamic
program, the state space is of size $O(n)$. Moreover, the recursion
formula requires going through each edge once, which is also of size
$O(n)$ because each product can only transit to at most one other
product. Therefore, the computational complexity of the DP is
$O(n)$. $\hfill\Box$\\

{\noindent\bf Proof of Theorem \ref{th:MCST_app}.}
Let $S^*$ be the optimal assortment under the MCST model with
optimal value $\OPT$ and the best recommend sets $R^{S^*}_j$ for all
$j\notin S^*$. Let $S_t =\{1,...,t\}$ be the optimal revenue-ordered
assortment with revenue $\mathbf{R}^{ro}$ and the best recommended
sets $R^{S_t}_j$ for $j > t$. First, we have
\begin{alignat}{2}
\OPT = & \sum_{i = 1}^{n}r_i\mathrm{Pr}_i(S^*, R^{S^*}_j) =
\sum^{n}_{t = 1}(r_{t} - r_{t+1})\sum_{i = 1}^{t}\mathrm{Pr}_i(S^*,
R^{S^*}_j) \nonumber\\
&= \sum^{n}_{t = 1}(r_{t} - r_{t+1})\sum_{i\in S^*\cap
S_t}\mathrm{Pr}_i(S^*, R^{S^*}_j)
 = \sum^{n}_{t = 1}\frac{r_{t} - r_{t+1}}{r_t}\sum_{i\in S^*\cap S_t}r_t\mathrm{Pr}_i(S^*, R^{S^*}_j). \label{eq:MCST_roapp_optlog}
\end{alignat}
Note that we let $r_{n+1} = 0$ and the equalities hold since
$\mathrm{Pr}_i(S^*, R^{S^*}_j) = 0$ if $i\not\in S^*$. From
\eqref{eq:MCST_choice_rj}, we have
\begin{alignat}{2}
&  \sum_{i\in S^*\cap S_t}r_t\mathrm{Pr}_i(S^*, R^{S^*}_j) \nonumber\\
& =  \sum_{i\in S^*\cap S_t}\lambda_ir_t + \sum_{i\in S^*\cap S_t}\sum_{j:j\not\in S, i\in R^{S^*}_j}\lambda_jr_t\frac{v_{ji}}{\sum_{i\in R^{S^*}_j}v_{ji} + v_{j0}} \nonumber\\
& = \sum_{i\in S^*\cap S_t}\lambda_ir_t + \sum_{j:j\not\in S^*}\lambda_jr_t\frac{\sum_{i\in R^{S^*}_j\cap S_t}v_{ji}}{\sum_{i\in R^{S^*}_j}v_{ji} + v_{j0}} \nonumber\\
& \leq \sum_{i\in S^*\cap S_t}\lambda_ir_t + \sum_{j:j\not\in S^*\cap S_t}\lambda_jr_t\frac{\sum_{i\in R^{S^*}_j\cap S_t}v_{ji}}{\sum_{i\in R^{S^*}_j\cap S_t}v_{ji} + v_{j0}}\label{eq:MCST_roapp_opt2_1}\\
& = \sum_{i\in S^*\cap S_t}\lambda_ir_t + \sum_{j:j\not\in S_t}\lambda_jr_t\frac{\sum_{i\in R^{S^*}_j\cap S_t}v_{ji}}{\sum_{i\in R^{S^*}_j\cap S_t}v_{ji} + v_{j0}} + \sum_{j:j\in S_t\setminus S^*}\lambda_jr_t\frac{\sum_{i\in R^{S^*}_j\cap S_t}v_{ji}}{\sum_{i\in R^{S^*}_j\cap S_t}v_{ji} + v_{j0}}\nonumber\\
& \leq \sum_{i\in S^*\cap S_t}\lambda_ir_i + \sum_{j:j\not\in S_t}\lambda_j\frac{\sum_{i\in R^{S^*}_j\cap S_t}v_{ji}r_i}{\sum_{i\in R^{S^*}_j\cap S_t}v_{ji} + v_{j0}} + \sum_{j:j\in S_t\setminus S^*}\lambda_jr_j\frac{\sum_{i\in R^{S^*}_j\cap S_t}v_{ji}}{\sum_{i\in R^{S^*}_j\cap S_t}v_{ji} + v_{j0}} \label{eq:MCST_roapp_opt2_2}\\
& \leq \sum_{i\in S^*\cap S_t}\lambda_ir_i + \sum_{j:j\not\in S_t}\lambda_j\frac{\sum_{i\in R^{S_t}_j}v_{ji}r_i}{\sum_{i\in R^{S_t}_j}v_{ji} + v_{j0}} + \sum_{j:j\in S_t\setminus S^*}\lambda_jr_j  \label{eq:MCST_roapp_opt2_3}\\
&= \sum_{i\in S_t}\lambda_ir_i + \sum_{j:j\not\in
S_t}\lambda_j\frac{\sum_{i\in R^{S_t}_j}v_{ji}r_i}{\sum_{i\in
R^{S_t}_j}v_{ji} + v_{j0}} = \sum_{i\in S_t}r_i\mathrm{Pr}_i(S_t,
R^{S_t}_j) = \mathbf{R}^{ro}, \label{eq:MCST_roapp_opt2}
\end{alignat}
where \eqref{eq:MCST_roapp_opt2_1} holds since $S^*\cap S_t\subseteq
S^*$ and $R^{S^*}_j\cap S_t\subseteq R^{S^*}_j$,
\eqref{eq:MCST_roapp_opt2_2} holds since $r_t \leq r_i$ for each
$i\in S_t$, \eqref{eq:MCST_roapp_opt2_3} holds since
$\frac{\sum_{i\in R^{S^*}_j\cap S_t}v_{ji}}{\sum_{i\in R^{S^*}_j\cap
S_t}v_{ji} + v_{j0}}\leq 1$, and the set $R^{S_t}_j$ is an optimal
solution to $\displaystyle\max\limits_{R_j\subseteq
S_t}\frac{\sum_{i\in R_j}r_iv_{ji}}{\sum_{i\in R_j}v_{ji} + v_{j0}}$
by definition and $R^{S^*}_j\cap S_t \subseteq S_t$. Note that
\eqref{eq:MCST_roapp_opt2} is the revenue of the assortment $S_t$
under its best recommended set. Plugging it into
\eqref{eq:MCST_roapp_optlog}, and the revenue-ordered assortment
returns the assortment with maximum revenue among all $S_t$, we
obtain \begin{equation} \OPT \leq \sum^{n}_{t = 1}\frac{r_{t} -
r_{t+1}}{r_t}\mathbf{R}^{ro} \leq \left(1 + \sum_{t =
1}^{n-1}\int^{r_t}_{r_{t+1}}\frac{1}{x} dx\right) \mathbf{R}^{ro} =
\left(1+\log \frac{r_{1}}{r_{n}}\right)\mathbf{R}^{ro}.
\label{eq:MCST_roapp_optlog2} \end{equation} Furthermore, let $\tilde{r}_1, ..., \tilde{r}_d$
be the distinct revenue values sorted in decreasing order.
Similarly, we can obtain from \eqref{eq:MCST_roapp_opt2} that
{\small\begin{equation} \OPT = \sum_{i \in S^*}r_i\mathrm{Pr}_i(S^*, R^{S^*}_j) =
\sum_{t = 1}^{d}\tilde{r}_t\sum_{i \in S^*, r_{i} = \tilde{r}_{t}}\mathrm{Pr}_i(S^*,
R^{S^*}_j) \leq \sum_{t = 1}^{d}\tilde{r}_t\sum_{i \in S^* \cap
S_t}\mathrm{Pr}_i(S^*, R^{S^*}_j) \leq d\mathbf{R}^{ro}.
\label{eq:MCST_roapp_optk} \end{equation}} Combining
\eqref{eq:MCST_roapp_optlog2} and \eqref{eq:MCST_roapp_optk}, we
obtain that the revenue of the best revenue-ordered assortment
achieves at least $\max\left\{\frac{1}{d},\frac{1}{1+\log
\frac{r_{\max}}{r_{\min}}}\right\}$ fraction of the optimal revenue.

\begin{figure}[htbp]
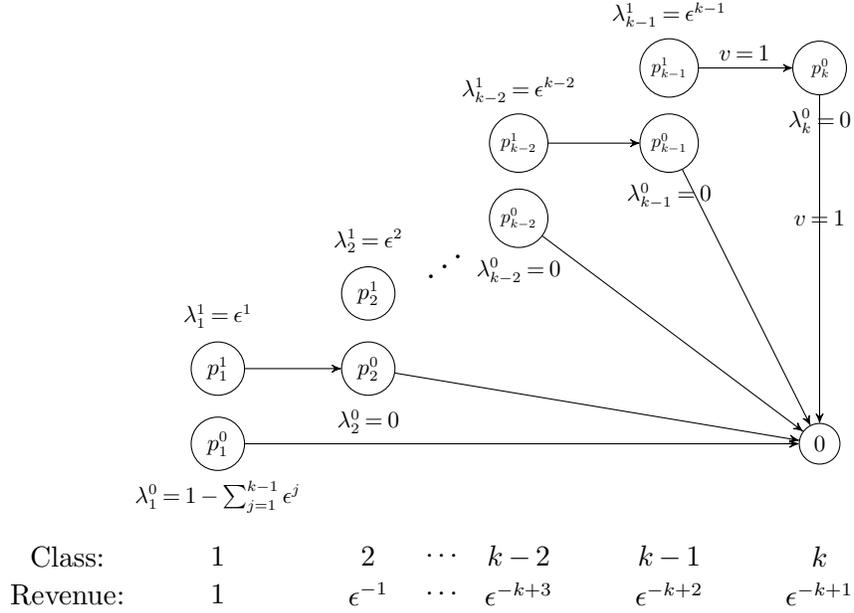

\begin{center}
\tikz{ \node(no) at (10, 0) [circle,scale=.8,draw] {$0$}; \node(pk0)
at (10, 5) [circle,scale=.65,draw] {$\;\;p^0_{k\;\;}$}
    edge[->,>=stealth'] (no);
\node(pkm11) at (8, 5) [circle,scale=.65,draw] {$p^1_{k-1}$}
    edge[->,>=stealth'] (pk0);
\node(pkm10) at (8, 4) [circle,scale=.65,draw] {$p^0_{k-1}$}
    edge[->,>=stealth'] (no);
\node(pkm21) at (6, 4) [circle,scale=.65,draw] {$p^1_{k-2}$}
    edge[->,>=stealth'] (pkm10);
\node(pkm20) at (6, 3) [circle,scale=.65,draw] {$p^0_{k-2}$}
    edge[->,>=stealth'] (no);
\node at (5,2.5) [scale=1.2] {$\iddots$}; \node(p21) at (4,2)
[circle,scale=.8,draw] {$p^1_{2}$}; \node(p20) at (4, 1)
[circle,scale=.8,draw] {$p^0_{2}$}
    edge[->,>=stealth'] (no);
\node(p11) at (2, 1) [circle,scale=.8,draw] {$p^1_{1}$}
    edge[->,>=stealth'] (p20);
\node(p10) at (2, 0) [circle,scale=.8,draw] {$p^0_{1}$}
    edge[->,>=stealth'] (no);
\node at (0, -1.5) {Class:};\node at (2, -1.5) {$1$};\node at (4,
-1.5) {$2$};\node at (6, -1.5) {$k-2$}; \node at (8, -1.5)
{$k-1$};\node at (10, -1.5) {$k$}; \node at (5, -1.5) {$\cdots$};
\node at (0, -2) {Revenue:};\node at (2, -2) {$1$};\node at (4, -2)
{$\epsilon^{-1}$};\node at (6, -2) {$\epsilon^{-k+3}$};\node at (8,
-2) {$\epsilon^{-k+2}$};\node at (10, -2) {$\epsilon^{-k+1}$}; \node
at (5, -2) {$\cdots$}; \node at (10, 4+.3) [scale=.8]{$\lambda^0_{k}
= 0$}; \node at (9, 5.2) [scale=.8]{$v = 1$}; \node at (10, 3)
[scale=.8]{$v = 1$}; \node at (8, 6-.3) [scale=.8]{$\lambda^1_{k-1}
= \epsilon^{k-1}$};\node at (8, 3+.3) [scale=.8]{$\lambda^0_{k-1} =
0$}; \node at (6, 5-.3) [scale=.8]{$\lambda^1_{k-2} =
\epsilon^{k-2}$};\node at (6, 2+.3) [scale=.8]{$\lambda^0_{k-2} =
0$}; \node at (4, 3-.3) [scale=.8]{$\lambda^1_{2} =
\epsilon^{2}$};\node at (4, 0+.3) [scale=.8]{$\lambda^0_{2} = 0$};
\node at (2, 2-.3) [scale=.8]{$\lambda^1_{1} = \epsilon^{1}$};\node
at (2, -1+.3) [scale=.8]{$\lambda^0_{1} = 1 - \sum_{j =
1}^{k-1}\epsilon^j$}; }
\end{center}
\caption{Tight Example of Revenue-Ordered Assortment}
\label{fig:MCST_tight}
\end{figure}

Finally, we show that the performance bound is tight in the sense
that we can find an instance such that the revenue-ordered assortment
only achieves $O\left(\frac{1}{d}\right)$ and
$O\left(\frac{1}{1+\log \frac{r_{\max}}{r_{\min}}}\right)$ fraction
of the optimal revenue respectively. Let $k \geq 1$, $\epsilon \in
(0,1)$ and the number of products be $n = 2k - 1$. For ease of
exposition, we denote the products as $k$ classes, where each class
has two products except for the $k$th class which has only one product.
The revenue of product in class $i$ is $\epsilon^{-i+1}$, for $i =
1, ..., k$. For $i = 1, ..., k-1$, one product in class $i$ is
denoted as product $p^1_i$, which has arrival probability
$\epsilon^i$; and the other product is denoted as $p^0_i$, each has
arrival probability $0$, except for the product $p^0_1$ which has
arrival probability $1 - \sum_{j = 1}^{k-1}\epsilon^j$. The class
$k$ has only product $p^0_k$, which has arrival probability 0. For
product $p^1_i$, $i = 1, ..., k-1$, the transition probability to
the product $p^0_{i+1}$ is 1, and the transition probabilities to
other products are 0 (including to the no-purchase option). For
product $p^0_i$, $i = 1, ..., k$, the transition probability to the
no-purchase option is 1, and the transition probabilities to other
products are 0. In other words, except for the transitions to the
no-purchase option, the only possible transition is from an
unavailable product $p^1_i$ to an available product $p^0_{i+1}$. The
construction of the case is illustrated in Figure
\ref{fig:MCST_tight}.

We first discuss the bound $O(1/d)$. Note that this instance has
exactly $k$ distinct revenue values, and within each class, there
are three possibilities to choose the products. Let $S_t$ be the
revenue-ordered assortment that all the products in classes $k$,
$k-1$, ..., $k-t+1$ are selected, $\forall t = 1,..., k$, and $S^0_t
= S_t\setminus \{p^1_{k-t+1}\}$ and $S^1_t = S_t \setminus
\{p^0_{k-t+1}\}$ be the other two types of revenue-ordered
assortments, $\forall t = 2,..., k$. It can be computed that for $t
= 1, ..., k-1$, the revenue of assortment $S_t$ is $\sum^{k-1}_{j =
k - t+1}\epsilon^{j}\cdot\epsilon^{-j+1} +
\epsilon^{k-t}\cdot\epsilon^{-k+t} = 1 + (t - 1)\epsilon$, where the
first term in the left-hand side corresponds to the $p^1_j$ products
for $j =  k-t+1, ..., k-1$ selected in $S_t$, and the second term
corresponds to the unavailable product $p^1_{k-t}$ transiting to the
product $p^0_{k - t+1}$ in $S_t$ (the revenue for selecting $p^0_j$
products is 0 since the arrival probability is 0); the revenue of
the last revenue-ordered assortment $S_k$ is $\sum^{k-1}_{j =
1}\epsilon^{j}\cdot\epsilon^{-j+1} + (1 - \sum_{j =
1}^{k-1}\epsilon^j)\cdot \epsilon^{-1+1} = 1 + (k - 2)\epsilon -
\sum_{j = 2}^{k-1}\epsilon^j$. For $t = 2,..., k$, the revenue of
assortment $S^1_t$ is $\sum^{k-1}_{j = k -
t+1}\epsilon^{j}\cdot\epsilon^{-j+1} = (t - 1)\epsilon$, and the
revenue of assortment $S^0_t$ is $\sum^{k-1}_{j = k -
t+2}\epsilon^{j}\cdot\epsilon^{-j+1} +
\epsilon^{k-t+1}\cdot\epsilon^{-k+t+1} +
\epsilon^{k-t}\cdot\epsilon^{-k+t} = 2 + (t - 2)\epsilon$, except
that $S^0_k = 2 + (k - 2)\epsilon - \sum_{j = 2}^{k-1}\epsilon^j$.
If we select all the $p^0_i$ products for $i = 1, ..., k$ but not
the $p^1_i$ products, then the revenue is $\sum^{k}_{j =
2}\epsilon^{j-1}\cdot\epsilon^{-j+1} + (1 - \sum_{j =
1}^{k-1}\epsilon^j)\cdot \epsilon^{-1+1} = k -\sum_{j =
1}^{k-1}\epsilon^j$. Therefore, the ratio of the revenue of the
revenue-ordered assortment and the optimal assortment is $2/k +
O(\epsilon)$, it goes to $2/k = O\left(1/d\right)$ when $\epsilon$
goes to 0. Moreover, we can see that
$\log{(\frac{r_{\max}}{r_{\min}})}= (k - 1)\log{(1/\epsilon)}$.
Thus, as $k$ goes to infinity, the second bound is
also tight. $\hfill\Box$\\

{\noindent\bf Proof of Theorem \ref{th:mip_equiv}.} Given any fixed
$\bm{x}\in \{0, 1\}^n$, we prove that the remaining optimization
problem in \eqref{eq:MCST_mlp2} and \eqref{eq:MCST_mlp1} are the
same. First, fixing $\bm{x}$, we show that \eqref{eq:MCST_mlp1} is
equivalent to the problem
\begin{subequations}
\label{eq:MCST_mlp1r}
\begin{alignat}{2}
\max\;\;\; & \sum_{j=1}^n\lambda_j(1 - x_j)\sum^{n}_{i = 1}(r_i - z^*_j)\frac{v_{ji}}{v_{j0}}y_{ji} \label{eq:MCST_mlp1r_obj}\\
\mathrm{s.t.}\;\;\;& y_{ji} \leq x_{i},& \qquad\forall i, j = 1, \ldots , n, \label{eq:MCST_mlp1r_c1}\\
& y_{ji} \leq 1 - x_{j},& \qquad\forall i, j = 1, \ldots , n, \label{eq:MCST_mlp1r_c2}\\
& 0\leq y_{ji} \leq 1, & \qquad\forall i, j = 1, \ldots ,
n,\nonumber
\end{alignat}
\end{subequations}
where $\bm{y}^*$ is an optimal solution to \eqref{eq:MCST_mlp1}
(while fixing $\bm{x}$) and $z^*_{j} = \sum_{i =
1}^{n}{r_iv_{ji}y^*_{ji}}/(\sum_{i = 1}^{n}v_{ji}y^*_{ji} +
v_{j0})$. By the definition of $\bm{y}^*$ and $z^*_j$, we have
$\sum_{i = 1}^{n}r_iv_{ji}y^*_{ji} = z^*_j\left(\sum_{i =
1}^{n}v_{ji}y^*_{ji} + v_{j0}\right)$, that is, $\sum^n_{i = 1}(r_i
- z^*_j)\frac{v_{ji}}{v_{j0}}y^*_{ji} = z^*_{j}$. Therefore, the
optimal value to \eqref{eq:MCST_mlp1r} is at least as large as the
optimal value to \eqref{eq:MCST_mlp1}. Now let $\tilde{\bm{y}}$ be
an optimal solution to \eqref{eq:MCST_mlp1r}. Note that
\eqref{eq:MCST_mlp1r} is a linear programming problem. We claim that
the constraint matrix of $\bm{y}$ formed by \eqref{eq:MCST_mlp1r_c1}
and \eqref{eq:MCST_mlp1r_c2} is a totally unimodular matrix with
dimension $2n^2 \times n^2$, since there is exactly one nonzero
entry in each row. Therefore, we can assume that $\tilde{\bm{y}}\in
\{0, 1\}^{n^2}$ (see Nemhauser and Wolsey 1988)\footnote{G. L.
Nemhauser and L. A. Wolsey. 1988. {\it Integer and Combinatorial
Optimization}. Wiley-Interscience, New York, NY, USA.}. Note that
for any $j$,
$$\lambda_j(1 - x_j)z^*_j = \lambda_j(1 - x_j)\frac{\sum_{i =
1}^{n}r_iv_{ji}y^*_{ji}}{\sum_{i = 1}^{n}v_{ji}y^*_{ji} + v_{j0}}
\geq \lambda_j(1 - x_j)\frac{\sum_{i =
1}^{n}r_iv_{ji}\tilde{y}_{ji}}{\sum_{i = 1}^{n}v_{ji}\tilde{y}_{ji}
+ v_{j0}}.$$ By rearranging terms, we get
$$\lambda_j(1 - x_j)z^*_j \geq \lambda_j(1 - x_j)\sum^{n}_{i = 1}(r_i
- z^*_j)\frac{v_{ji}}{v_{j0}}\tilde{y}_{ji},$$ thus $\sum_{j=1}^n
\lambda_j(1 - x_j)z^*_j \geq \sum_{j=1}^n \lambda_j(1 -
x_j)\sum^{n}_{i = 1}(r_i -
z^*_j)\frac{v_{ji}}{v_{j0}}\tilde{y}_{ji}$. Therefore, the optimal
value to \eqref{eq:MCST_mlp1} is at least as large as the optimal
value to \eqref{eq:MCST_mlp1r}, i.e., the optimal values to the two
problems are the same. Thus, for any given $\bm{x}$ and given an
optimal solution to one of the two problems, we can construct a
feasible solution to the other which has the same objective value,
which implies that \eqref{eq:MCST_mlp1} and \eqref{eq:MCST_mlp1r}
are equivalent.

Next we show that \eqref{eq:MCST_mlp1r} and \eqref{eq:MCST_mlp2} are
equivalent. On one hand, let $\bm{y}^*$ be an optimal solution to
\eqref{eq:MCST_mlp1r}, and $\OPT(\bm{y}^*)$ be the optimal value. We
construct the solution $z_{ji} = {v_{ji}y^*_{ji}(1 - x_j)}/(\sum_{i
= 1}^{n}v_{ji}y^*_{ji} + v_{j0})$ and $z_{j0} = {v_{j0}(1 -
x_j)}/(\sum_{i = 1}^{n}v_{ji}y^*_{ji} + v_{j0})$, $\forall i, j = 1,
..., n$ to problem \eqref{eq:MCST_mlp2}. It can be verified that the
solution $\bm{z}$ is feasible to \eqref{eq:MCST_mlp2}. In
particular, if $x_i = 0$, then $y^*_{ji} = 0$ from constraint
\eqref{eq:MCST_mlp1r_c1} and we have $z_{ji} = 0 \leq x_i$;
otherwise, by the above definition and $x_j \in \{0, 1\}$, $z_{ji}
\leq 1 - x_j \leq 1 = x_i$, therefore constraint
\eqref{eq:MCST_mlp2_c1} is always satisfied. Constraint
\eqref{eq:MCST_mlp2_c3} is clearly satisfied from the definition of
$\bm{z}$. Moreover, we have $\displaystyle\frac{z_{ji}}{z_{j0}} =
\frac{v_{ji}y^*_{ji}}{v_{j0}}$, i.e., $\displaystyle z_{ji} =
\frac{v_{ji}}{v_{j0}}z_{j0}y^*_{ji}$, it follows that $\displaystyle
z_{ji} \leq \frac{v_{ji}}{v_{j0}}z_{j0}$ since $y^*_{ji} \leq 1$ and
thus constraint \eqref{eq:MCST_mlp2_c4} is satisfied. Note that for
each $j = 1, \ldots, n$, $\sum_{i = 1}^{n}r_iz_{ji} = (1 -
x_j)\sum_{i = 1}^{n}r_iv_{ji}y^*_{ji}/(\sum_{i =
1}^{n}v_{ji}y^*_{ji} + v_{j0})$, Thus it implies that the objective
value of $\bm{z}$ to \eqref{eq:MCST_mlp2} is at least
$\OPT(\bm{y}^*)$. Therefore, the optimal value to
\eqref{eq:MCST_mlp2} is at least as large as the optimal value to
\eqref{eq:MCST_mlp1r}.

On the other hand, let $\tilde{\bm{z}}$ be an optimal solution to
\eqref{eq:MCST_mlp2}, and $\OPT(\tilde{\bm{z}})$ be the optimal
value. We construct the solution $y_{ji} =
v_{j0}\tilde{z}_{ji}/(v_{ji}\tilde{z}_{j0})$, $\forall i, j = 1,
..., n$. Now we discuss its feasibility to \eqref{eq:MCST_mlp1r}. By
constraint \eqref{eq:MCST_mlp2_c4} and the construction of $y_{ji}$,
we have $0 \leq y_{ji} \leq 1$. For each $\tilde{z}_{ji}$, if
$\tilde{z}_{ji} = 0$, then $y_{ji} = 0$, and constraints
\eqref{eq:MCST_mlp1r_c1} and \eqref{eq:MCST_mlp1r_c2} must be
satisfied since $\bm{x} \in \{0, 1\}$; otherwise, $\tilde{z}_{ji} >
0$ and it implies that $x_i = 1$ and $1 - x_j = 1$, thus constraints
\eqref{eq:MCST_mlp1r_c1} and \eqref{eq:MCST_mlp1r_c2} are also
satisfied since $y_{ji} \leq 1$ as mentioned before. Note that in
the objective function \eqref{eq:MCST_mlp1r_obj}, the value $z^*_{j}
\geq \sum_{i = 1}^{n}{r_iv_{ji}y_{ji}}/(\sum_{i = 1}^{n}v_{ji}y_{ji}
+ v_{j0})$ by the definition of $z^*_{j}$. By rearranging terms, we
get $v_{j0}z^*_{j} \geq \sum^n_{i=1}r_iv_{ji}y_{ji} -
z^*_{j}\sum^n_{i=1}v_{ji}y_{ji} =
\sum^n_{i=1}r_iv_{j0}\tilde{z}_{ji}/\tilde{z}_{j0} -
z^*_{j}\sum^n_{i=1}v_{j0}\tilde{z}_{ji}/\tilde{z}_{j0}$, or
equivalently, $z^*_{j} \geq
\sum^n_{i=1}r_i\tilde{z}_{ji}/\tilde{z}_{j0} -
z^*_{j}\sum^n_{i=1}\tilde{z}_{ji}/\tilde{z}_{j0}$. Since
\eqref{eq:MCST_mlp1} and \eqref{eq:MCST_mlp1r} are equivalent and
$\OPT(\tilde{\bm{z}}) \geq \OPT(\bm{y}^*) = \sum_{j =
1}^{n}\lambda_jr_jx_j + \sum_{j=1}^n\lambda_j(1 - x_j)z^*_j$, the
optimal solution to \eqref{eq:MCST_mlp1r} is at least
\begin{eqnarray*}
\sum_{j = 1}^{n}\lambda_jr_jx_j + \sum_{j=1}^n\lambda_j(1 -
x_j)z^*_j \hspace{-.2cm}& \geq & \hspace{-.2cm}\sum_{j = 1}^{n}\lambda_jr_jx_j +
\sum_{j=1}^n\lambda_j(1 -
x_j)\left(\sum^n_{i=1}r_i\frac{\tilde{z}_{ji}}{\tilde{z}_{j0}} -
z^*_{j}\sum^n_{i=1}\frac{\tilde{z}_{ji}}{\tilde{z}_{j0}}\right)\\
& \geq & \hspace{-.2cm}\sum_{j = 1}^{n}\lambda_jr_jx_j + \sum_{j=1}^n\lambda_j(1 -
x_j)\sum^n_{i=1}r_i\frac{\tilde{z}_{ji}}{\tilde{z}_{j0}} - \lambda_j
\sum^{n}_{i =
1}r_i\tilde{z}_{ji}\sum^n_{i=1}\frac{\tilde{z}_{ji}}{\tilde{z}_{j0}}\\
& = &\hspace{-.2cm}\sum_{j = 1}^{n}\lambda_jr_jx_j + \sum_{j=1}^n\lambda_j(1 -
x_j)\sum^n_{i=1}\tilde{z}_{ji} = \OPT(\tilde{\bm{z}}),
\end{eqnarray*}
where the second last equality follows from constraint
\eqref{eq:MCST_mlp2_c3}. Therefore, the optimal value to
\eqref{eq:MCST_mlp1r} is at least as large as the optimal value to
\eqref{eq:MCST_mlp2}, and thus the two problems are equivalent.
Recall that the above arguments hold for arbitrary $\bm{x}\in \{0,
1\}^n$, thus the theorem holds. $\hfill\Box$\\

\end{document}